\documentclass[12pt, dvipsnames]{article}

\usepackage{authblk}
\usepackage{amssymb, amsmath, amsthm}
\usepackage{tikz}
\usepackage{mathtools}
\usepackage{stmaryrd}
\usepackage{booktabs}
\usepackage{import}
\usepackage{siunitx}
\usepackage[subpreambles=true]{standalone}
\usepackage{subcaption}
\usepackage{wrapfig}
\usepackage{bm}
\usepackage[numbers]{natbib}

\newtheorem{definition}{Definition}[section]
\newtheorem{lemma}[definition]{Lemma}

\newtheorem{remark}[definition]{Remark}
\newtheorem{assumption}[definition]{Assumption}

\newcommand{\diver}{\operatorname{div}}
\newcommand{\llcurly}{\{\!\{}
\newcommand{\rrcurly}{\}\!\}}
\newcommand{\jump}[1]{\llbracket #1 \rrbracket}
\newcommand{\mean}[1]{\llcurly #1 \rrcurly}
\newcommand{\dumux}{DuMu$^\text{x}$ }
\newcommand{\dunemmesh}{\texttt{dune-mmesh}}

\newcommand{\filliation}[5]{\affil{#2\vskip 0pt #3\vskip 0pt #4\vskip 0pt #5}}

\title{A Finite-Volume Moving-Mesh Method for Two-phase Flow in Fracturing Porous Media}

\author{Samuel Burbulla\thanks{Acknowledgement: Funded by the Deutsche Forschungsgemeinschaft (DFG, German Research Foundation) -- Project Number 327154368 -- SFB 1313.}
}
\author{Christian Rohde}

\date{}
\filliation{University of Stuttgart}{Institute of Applied Analysis and Numerical Simulation}{Pfaffenwaldring 57, 70569 Stuttgart, Germany}{samuel.burbulla@mathematik.uni-stuttgart.de}

\begin{document}

\maketitle

\begin{abstract}
Flow in fractured porous media is modeled frequently by discrete fracture-matrix approaches where fractures are treated as dimensionally reduced manifolds.
Generalizing earlier work we focus on two-phase flow in time-dependent fracture geometries including the fracture's aperture.
We present the derivation of a reduced model for immiscible two-phase flow in porous media. 
For the reduced model we present a fully conforming finite-volume discretization coupled with a moving-mesh method.
This method permits arbitrary movement of facets of the triangulation while being fully conservative.
In numerical examples we show the performance of the scheme and investigate the modeling error of the reduced model.

\end{abstract}

\section{Introduction}

Discrete fracture-matrix models are a widely-spread approach to model flow in fractured porous media. Hereby, fractures are reduced to lower-dimensional manifolds. On the one hand, this approach leads to a reduced computational effort and less geometrical issues especially when treating with very thin fractures. On the other hand, it can be useful to have an explicit representation of the fracture geometry. Especially in the context of propagating fractures, the fracture geometry becomes an unknown of the problem and a numerical discretization must be able to incorporate the movement.

In the majority of cases the fractures are represented geometrically by a lower-dimensional grid. Here, some approaches use non-conforming representations \cite{Fumagalli2013}, others follow a conforming ansatz where the fractures have to coincide with facets of the matrix mesh \cite{glaeser}. The non-conforming discretizations have the advantage to be independent of the surrounding grid geometry, but the fracture dynamics are more complicated to realize. On the other hand, in conforming methods, the geometries are coupled directly and fracture dynamics are almost trivial to realize. However, any movement results in re-meshing that has to be performed simultaneously on both grids to keep the interface conforming.

Besides that, we have to mention approaches with phase-field representation of the fractures \cite{sanghyun2018}. In contrast to discrete fracture models, the phase-field approaches are entirely flexible in geometry and do not require any re-meshing apart from adaptivity. However, they are computationally very expensive because a fine grid resolution is needed to resolve the phase-field gradient. In addition, no explicit geometrical representation of the fracture is given and all the complexity is shifted into the model.

In this work, we want to investigate a fully conforming discrete fracture approach. We believe that the additional effort in conforming re-meshing of the grid pays off in the easy coupling to the lower-dimensional grid. We derive a reduced model for capillarity-free two-phase flow in porous media and present a numerical scheme on the basis of a finite-volume discretization. The method is fully conservative, consistent for isotropic permeability tensors and does not require the construction of any cell stencils that might be expensive during re-meshing. The method also allows for a treatment of fully-resolved propagating fractures that will be used for the numerical analysis of the reduced model.
It is a new approach combining a mixed-dimensional modeling approach for two-phase flow in porous media with a finite-volume moving-mesh method for fracture propagation.

The paper is structured as follows. In Section 2, we describe the geometrical setting of our model and derive a reduced model for capillarity-free two-phase flow in porous media. In Section 3, we present the idea behind the finite-volume moving-mesh method together with the adjusted two-point flux discretization. In Section 4, we show the performance of the reduced model by presenting some numerical experiments where we compare the solutions of the reduced and the full dimensional model. In the end, we will summarize and discuss our results and give some outlook on future work.

\section{The Mathematical Model for Two-phase Flow in Fracturing Porous Media}

In this section we derive a mathematical model that governs capillarity-free two-phase flow in dynamically fracturing porous media.
As our primary approach we rely on the fractional flow formulation, see system (\ref{eq:fffg}) below, that is derived from the classical coupled formulation as e.g.\ in \cite{helmig1997multiphase}.

\subsection{Coupled Formulation}
\label{subsec:coupledform}

Let $\Omega \subset \mathbb{R}^n$ be open and bounded and $T > 0$. As described in \cite{helmig1997multiphase}, the dynamics of two incompressible and immiscible fluids (denoted as wetting $_w$ and non-wetting $_{nw}$ one) in porous media can be described by a system of differential equations, consisting of mass conservation laws and the subphases' Darcy laws, i.e.,
\begin{align}
 \label{eq:tps}
 \left.
 \begin{gathered}
 ( \phi \rho_\alpha S_\alpha )_t + \diver( \rho_\alpha v_\alpha ) = \rho_\alpha q_\alpha, \\
 \mathbf{v}_\alpha = - \lambda_\alpha(S_\alpha) \textbf{K} ( \nabla P_\alpha - \rho_\alpha \mathbf{g} )
 \end{gathered}
 \qquad \right. \quad
 \begin{gathered}
  \text{in } \Omega \times (0,T), \quad
  \alpha \in \{w, nw\}.
 \end{gathered}
\end{align}
Here, $S_\alpha: \Omega \times (0,T) \to [0,1]$ is the saturation of fluid $\alpha \in \{w, nw\}$, $P_\alpha: \Omega \times (0,T) \to \mathbb{R}$ is the corresponding phase pressure and $\mathbf{v}_\alpha: \Omega \times (0,T) \to \mathbb{R}^n$ stands for the fluid velocity.

The system (\ref{eq:tps}) is closed by the identities
\begin{subequations}
 \begin{align}
 \label{eq:tpsidentA}
  S_w + S_{nw} &= 1, \\
 \label{eq:tpsidentB}
  P_{nw} - P_w &= p_c(S_w),
 \end{align}
\end{subequations}
where $p_c$ is the capillary pressure function. Whereas (\ref{eq:tpsidentA}) is a natural saturation requirement, the equation (\ref{eq:tpsidentB}) is a widely-discussed assumption. For constitutive choices of $p_c$ depending on $S_w$ see \cite{helmig1997multiphase}. Actually, we will set $p_c \equiv 0$. This capillarity-free approach results in a hyperbolic-elliptic system of equations, see (\ref{eq:fffg}) below.
The given parameters in (\ref{eq:tps}), (\ref{eq:tpsidentA}), (\ref{eq:tpsidentB}) are the porosity $\phi = \phi(\mathbf{x}) \in (0,1]$, the constant phase density $\rho_\alpha \in \mathbb{R}^{>0}$, the symmetric and positive-definite intrinsic permeability tensor $\textbf{K} = \mathbf{K}(\mathbf{x})\in \mathbb{R}^{n \times n}$ and the gravitational acceleration vector $\mathbf{g} \in \mathbb{R}^n$.
The phase mobility function is given by $\lambda_\alpha( S_\alpha ) = k_\alpha( S_\alpha ) / \mu_\alpha$, where $k_\alpha = k_\alpha( S_\alpha )$ is the relative permeability (see Fig.\ \ref{fig:relperm}) and $\mu_\alpha \in \mathbb{R}$ is the dynamic viscosity of phase $\alpha$.
The function $q_\alpha: \Omega \times (0,T) \to \mathbb{R}$ is a source or sink term. We will specify the values for all these quantities when needed for the numerical experiments in Section \ref{sec:experiments}.

As mentioned above we will deviate from (\ref{eq:tps}) and rely our work on another equivalent formulation.

\subsection{The Fractional Flow Formulation}
\label{subsec:fractionalform}
Let us select the wetting phase saturation $S=S_w$ eliminating $S_{nw}$ by using (\ref{eq:tpsidentA}).
Introducing for $S$ the fractional flow function
\[ f(S) := \frac{\lambda_w(S)}{\lambda(S)}, \]
with $\lambda(S) := \lambda_w(S) + \lambda_{nw}(S)$ being the total mobility,
we reformulate the two-phase flow system (\ref{eq:tps}), (\ref{eq:tpsidentA}), (\ref{eq:tpsidentB}). Assuming zero capillary pressure, i.e. $p_c \equiv 0$, the two-phase flow system (\ref{eq:tps}) is equivalent to
\begin{align}
 \label{eq:fffg}
 \left.
 \begin{aligned}
 (\phi S)_t + \diver \mathbf{F}(S, \mathbf{v}) &= q_w, \\
 \mathbf{v} + \lambda(S) \mathbf{K} ( \nabla P - G(S) \mathbf{g} ) &= \mathbf{0}, \\
 \diver(\mathbf{v}) &= q_w + q_{nw}
 \end{aligned}
 \qquad \right. \quad
 \begin{gathered}
  \text{in } \Omega \times (0,T).
 \end{gathered}
\end{align}
In (\ref{eq:fffg}), the flux function $\mathbf{F} = \mathbf{F}(S, \mathbf{v})$ is defined as
\begin{equation*}
 \mathbf{F}(S, \mathbf{v}) := f(S) \mathbf{v} - f(S) \lambda_{nw}(S) \mathbf{K} ( \rho_{nw} - \rho_w ) \mathbf{g}.
\end{equation*}
For standard choices, $\mathbf{F}(\cdot, \mathbf{v})$ is non-monotone (see Fig.\ \ref{fig:fracflow}). The term $G$ is given by
\begin{equation*}
 G(S) := \frac{ \lambda_w(S) \rho_w + \lambda_{nw}(S) \rho_{nw} }{ \lambda(S) }.
\end{equation*}

In (\ref{eq:fffg}), the unknowns are the (wetting phase) saturation $S = S_w$, the total velocity $\mathbf{v} = \mathbf{v}_w + \mathbf{v}_{nw}$ and the global pressure $P = P_{nw} = P_w$. Note that $S_{nw}$ is computable from (\ref{eq:tpsidentA}).

It remains to put initial conditions $S^0: \Omega \to [0,1]$ for the saturation by $S(\mathbf{x}, 0) = S^0(\mathbf{x}), \mathbf{x} \in \Omega,$ and appropriate boundary conditions for $S$, $\mathbf{v}$ and/or $P$. We note that (\ref{eq:fffg}) is of mixed hyperbolic-elliptic type.

\subsection{The Discrete-Fracture Model}

Turning to fractured porous media, it is a wide-spread approach to model sufficiently thin fractures in porous media as lower-dimensional manifolds. Note, however, that we implicitly assume that the original fracture's aperture is clearly separated from the pore scale. Then, it is justified that the fracture persists on the Darcy scale. For single-phase flow, such discrete fracture models have been derived by transversal averaging \cite{JAFFRE2011967, VincentMartin2004}. We suggest a mathematical model to describe two-phase flow in fractured porous media on the basis of a discrete-fracture network approach by averaging a full dimensional fracture-bulk medium.

For $t \in [0,T]$ let us consider a connected open set $\Omega_f(t) \subset \Omega$, representing the original fracture. We suppose that the domain $\Omega$ is partitioned according to $\Omega = \big(\bar \Omega_b(t) \cup \bar \Omega_f(t)\big)^\circ$ (see Fig. \ref{fig:sketchdomain} for a sketch of the geometry), where $\Omega_b(t)$ denotes the bulk porous medium with $\Omega_b(t) \cap \Omega_f(t) = \emptyset$.

\begin{figure}[ht]
\begin{minipage}{0.49\textwidth}
 \centering
  \begin{tikzpicture} 

  \draw[thick, fill=black!30!white](4.38,1.8) to[out=165,in=30, distance=40] (0.8,2.8)
  						to[out=-20,in=165, distance=40] (4.0,1.3);  
  				
  \draw[thick] (0,1) to[out=150,in=190, distance=40] (1,4.7) 
  					 to[out=10, in=150] (4,3.55)
  					 to[out=-30, in=10, distance=40] (3,0.8)
  					 to[out=190, in=-30] (0,1);
  					 
  \node at (3.8,0.7) {$\Omega$};
  \node at (0.95,4.1) {$\Omega_b(t)$};
  \node[rotate=-25] at (2.4,2.35) {$\Omega_f(t)$};
 \end{tikzpicture}
 \caption{Sketch of the domain $\Omega$ divided into a bulk porous medium $\Omega_b(t)$ and the fracture region $\Omega_f(t)$. \\ \ \\}
 \label{fig:sketchdomain}
\end{minipage}
\hfill
\begin{minipage}{0.49\textwidth}
 \centering
  \begin{tikzpicture}
 
  \draw[thick, dashed](4.38,1.8) to[out=165,in=30, distance=40] (0.8,2.8)
  						to[out=-20,in=165, distance=40] (4.0,1.3);  
  						
  \draw[thick] (0,1) to[out=150,in=190, distance=40] (1,4.7) 
  					 to[out=10, in=150] (4,3.55)
  					 to[out=-30, in=10, distance=40] (3,0.8)
  					 to[out=190, in=-30] (0,1);
  					 
  \draw[very thick](0.8,2.8) to[out=5, in=165, distance=40] (4.2,1.55);
  
  \draw[<->] (2.48,1.97) -- (2.8,2.53) node[anchor=south west]{$d(\mathbf{s}, t)$};
  
  \draw[thick, ->] (1.7,2.7) -- +(-0.13,-0.52) node[anchor=north, yshift=-1pt]{$\mathbf{n}$};
 
  \node at (4.85,1.7) {$\gamma^+$};
  \node at (4.45,1.2) {$\gamma^-$};
  \node at (0.5,3.0) {$\Gamma(t)$};

 \end{tikzpicture}
 \caption{In the dimension-reduction ansatz the fracture domain $\Omega_f(t)$ is replaced by a lower-dimensional interface $\Gamma(t)$, together with an aperture $d = d(\mathbf{s}, t)$.}
 \label{fig:sketchdomainfrac}
\end{minipage}
\end{figure}

As in Sections \ref{subsec:coupledform}, \ref{subsec:fractionalform} we assume to have incompressible and immiscible two-phase flow in the porous medium for both the bulk domain $\Omega_b(t)$ and the fracture domain $\Omega_f(t)$. However, the porosities $\phi, \phi^f$, the intrinsic permeabilities $\mathbf{K}, \mathbf{K}^f$ and the relative permeabilities $k_\alpha, k_\alpha^f$ might differ.

\begin{remark}
Within fractures a porous media (two-phase) flow occurs if they are filled with debris. Effective porous media flow can also be induced by wall roughnesses. We assume that this behaviour can be characterized by the physical properties of a porous medium in the fracture.
\end{remark}

To formulate a closed model it remains to impose suitable coupling conditions at the interfaces between the subdomains. For this purpose, define the space-time sets
\begin{align}
 \Omega_b^t &\coloneqq \{ \Omega_b(t) \times \{t\} \mid t \in [0,T] \},
 \\
 \Omega_f^t &\coloneqq \{ \Omega_f(t) \times \{t\} \mid t \in [0,T] \}.
\end{align}

Adapting the fractional-flow formulation from (\ref{eq:fffg}), the fluid states \\ $S^i: \Omega_i^t \to [0,1]$, $P^i: \Omega_i^t \to \mathbb{R}$ and $\mathbf{v}^i: \Omega_i^t \to \mathbb{R}^n$ are governed by
\begin{align}
 \label{eq:fffgequidim}
 \left.
 \begin{aligned}
  (\phi^i S^i)_t + \diver \mathbf{F}^i(S^i, \mathbf{v}^i) &= q_w, \\
  \mathbf{v}_i + \lambda^i(S^i) \mathbf{K}^i (\nabla P^i - G^i(S^i) \mathbf{g}) &= 0, \\
  \diver(\mathbf{v}^i) &= q_w + q_{nw}
 \end{aligned}
 \qquad \right. \quad
 \begin{gathered}
  \text{in } \Omega_i^t,\quad  i \in \{b,f\},
 \end{gathered}
\end{align}
with the coupling conditions
\begin{subequations}
\begin{align}
 \label{eq:couplcondsA}
 \mathbf{F}^b(S^b, \mathbf{v}^b) \cdot \mathbf{n} &= \mathbf{F}^f(S^f, \mathbf{v}^f) \cdot \mathbf{n}, \\
 \label{eq:couplcondsB}
 P^b &= P^f, \\
 \label{eq:couplcondsC}
 \mathbf{v}^b \cdot \mathbf{n} &= \mathbf{v}^f \cdot \mathbf{n}
\end{align}
\end{subequations}
on $\partial \Omega_b(t) \cap \partial \Omega_f(t)$ for all $t \in (0,T)$ and $i \in \{b,f\}$.

Here, the physical quantities are expected to differ in the two subdomains and, therefore, are indicated by the superscript $i \in \{b, f\}$. Appropriate initial data and boundary conditions have to be added.

The coupling conditions (\ref{eq:couplcondsA})-(\ref{eq:couplcondsC}) are a natural choice such that the problem (\ref{eq:fffgequidim}) is equivalent to the problem (\ref{eq:fffg}), defined on the complete domain $\Omega = \Omega_b \cup \Omega_f$.

Now, let us replace the fracture $\Omega_f(t)$ by a centered $(n-1)$-dimensional hypersurface $\Gamma(t)$ and a given aperture function $d \in \mathcal{C}^1(\Gamma^t)$ with $d > 0$, where $\Gamma^t \coloneqq \{ \Gamma(t) \times \{t\} \mid t \in (0,T) \}$.

For this, we require a simple geometrical setup with no junctions to be present, where the fracture boundaries are given by the graph of the aperture function defined on $\Gamma(t)$. Let us assume that the fracture does not interact with the boundary, i.e., $\overline{\Omega_f(t)} \subset \Omega^\circ$.

We shall orient the two sides of the surface. To do so, we choose a normal $\mathbf{n} = \mathbf{n}(\mathbf{s}, t) \in \mathcal{S}^n$ in $\mathbf{s} \in \Gamma(t)$. The sides will be denoted by positive (+) and negative (-) signs where $\mathbf{n}$ is assumed to point to the negative side. Later, we will also use the notation $\mathbf{n} = \mathbf{n}^+ = -\mathbf{n}^-$.

Let us consider smooth solutions $S^i$, $P^i$ and $\mathbf{v}^i$ for $i \in \{ b,f \}$ of system (\ref{eq:fffgequidim}) satisfying the coupling conditions (\ref{eq:couplcondsA})-(\ref{eq:couplcondsC}). Then, the reduced model is obtained by averaging along the line segments
\begin{equation*}
  L(\mathbf{s},t) := \Bigg\{ \mathbf{s} + \frac{r}{2} \mathbf{n} \Biggm| r \in \big(-d(\mathbf{s}, t), d(\mathbf{s}, t)\big) \Bigg\} \subset \Omega^f(t), \qquad \mathbf{s} \in \Gamma(t).
\end{equation*}

We introduce the projection matrices $\mathbf{N} \coloneqq \mathbf{n} \otimes \mathbf{n} \in \mathbb{R}^{n \times n}$ and $\mathbf{T} \coloneqq \mathbf{I} - \mathbf{N} \in \mathbb{R}^{n \times n}$. For a differentiable function $\varphi: \Gamma \to \mathbb{R}$ we define tangential and normal derivations by
\begin{align}
 \nabla_\tau \varphi &\coloneqq \nabla ( \mathbf{T} \bar \varphi ), &
 \nabla_\mathbf{n} \varphi &\coloneqq \nabla ( \mathbf{N} \bar \varphi ), \\
 \diver_\tau \varphi &\coloneqq \diver ( \mathbf{T} \bar \varphi ), &
 \diver_\mathbf{n} \varphi &\coloneqq \diver ( \mathbf{N} \bar \varphi ),
\end{align}
where $\bar \varphi$ is an extension of $\varphi$ to an open set including $\Gamma$. To each $\mathbf{s} \in \Gamma(t)$ we associate the boundary points
\begin{align}
 \mathbf{s}_\pm(\mathbf{s}, t) \coloneqq \mathbf{s} \mp \frac{1}{2} d(\mathbf{s}, t) \mathbf{n}.
\end{align}
Denote by
\begin{align}
 \gamma_\pm(t) \coloneqq \Bigg\{ \mathbf{s}_\pm(\mathbf{s}, t) \Biggm| \mathbf{s} \in \Gamma \Bigg\} \quad
\end{align}
the boundary surfaces on each side of the fracture (see Fig.\ \ref{fig:sketchdomainfrac}).
Using this notation, we define for some function $\varphi: \Omega_b \to \mathbb{R}$ the traces
\begin{align}
 \varphi^\pm = \varphi^\pm(\mathbf{s}) \coloneqq \lim_{\substack{\varepsilon \to 0 \\ \varepsilon > 0}} \varphi( \mathbf{s}_\pm \mp \varepsilon \mathbf{n} ),
\end{align}
and the jump and mean values by
\begin{align}
  \llbracket \varphi \rrbracket \coloneqq \varphi^+ - \varphi^- \qquad \text{ and }\qquad
  \llcurly \varphi \rrcurly \coloneqq \frac{\varphi^+ + \varphi^-}{2}.
\end{align}

We proceed and define reduced quantities on $\Gamma(t)$ by averaging along the line segments $L(\mathbf{s}, t)$, namely
\begin{align}
  \label{eq:avgquansA}
  \mathbf{v}_\Gamma(\mathbf{s}, t) \coloneqq \frac{1}{d(\mathbf{s}, t)} \int_{L(\mathbf{s}, t)} \mathbf{T} \mathbf{v}^f(\cdot, t)\ dS, \\
  \label{eq:avgquansB}
  S_\Gamma(\mathbf{s}, t) \coloneqq \frac{1}{d(\mathbf{s}, t)} \int_{L(\mathbf{s}, t)} S^f(\cdot, t)\ dS, \\
  \label{eq:avgquansC}
  P_\Gamma(\mathbf{s}, t) \coloneqq \frac{1}{d(\mathbf{s}, t)} \int_{L(\mathbf{s}, t)} P^f(\cdot, t)\ dS.
\end{align}

To simplify the model derivation, we assume that $\phi^f$ and $\mathbf{K}^f$ are constant along each $L(\mathbf{s}, t)$. For smooth functions $k_\alpha^f$, $f^f$ and $G^f$, that may depend non-linearly on $S^f$, we approximate all evaluation in $L(\mathbf{s}, t)$ by the evaluation at the mean quantity, i.e.,
\begin{align}
\label{eq:relativeperdmdefA}
k^f_\alpha(S^f(\mathbf{x}, t)) &= k^f_\alpha(S_\Gamma(\mathbf{s}, t)) + O(\bar d), \\
\label{eq:relativeperdmdefB}
f^f(S^f(\mathbf{x}, t)) &= f^f(S_\Gamma(\mathbf{s}, t)) + O(\bar d), \\
\label{eq:relativeperdmdefC}
G^f(S^f(\mathbf{x}, t)) &= G^f(S_\Gamma(\mathbf{s}, t)) + O(\bar d),
\end{align}
for all $\mathbf{x} \in L(\mathbf{s}, t)$ where $\bar d \coloneqq \max \{ d(\mathbf{s}, t) \mid t \in [0,T], \mathbf{s} \in \Gamma(t) \}$.

\begin{remark}
The affine choice for $k^f_\alpha$, i.e., $k^f_w(S^f) = S^f$ and $k^f_{nw}(S^f) = 1-S^f$, implies that conditions (\ref{eq:relativeperdmdefA})-(\ref{eq:relativeperdmdefC}) are satisfied exactly. The $O(\bar d)$-term vanishes.
\end{remark}

The further derivation of our model needs the following formula.
\begin{lemma}
For $\varphi \in \mathcal{C}^1(\Omega_f^t)$ we have for all $(\mathbf{s}, t) \in \Gamma^t$ the identity
\begin{align}
 \begin{split}
 \label{eq:timeintegral}
 &\frac{1}{d(\mathbf{s}, t)} \int_{L(\mathbf{s}, t)} \frac{d}{dt} \varphi(\mathbf{x}, t) \ dS \\
 &\qquad = \frac{d}{dt} \varphi_\Gamma(\mathbf{s}, t)
  + \frac{ \partial_t d(\mathbf{s}, t) }{ d(\mathbf{s}, t) } \Big( \varphi_\Gamma(\mathbf{s}, t) - \frac{\varphi \rvert_{\mathbf{s}_-} + \varphi \rvert_{\mathbf{s}_+}}{2} \Big),
 \end{split}
\end{align}
with $\varphi_\Gamma(\mathbf{s}, t) \coloneqq \frac{1}{d(\mathbf{s}, t)} \int_{L(\mathbf{s}, t)} \varphi(\mathbf{x}, t) \ dS$.
\begin{proof}
We compute
\begin{align*}
 &\frac{d}{dt} \varphi_\Gamma(\mathbf{s}, t) \\
 &\quad = \frac{d}{dt}  \Big(\frac{1}{d(\mathbf{s}, t)}\Big) \int_{L(\mathbf{s}, t)} \varphi(\mathbf{x}, t) \ dS + \frac{1}{d(\mathbf{s}, t)} \frac{d}{dt} \Big(\int_{L(\mathbf{s}, t)} \varphi(\mathbf{x}, t) \ dS\Big) \\
 &\quad = - \frac{\partial_t d(\mathbf{s}, t)}{d(\mathbf{s}, t)^2} \int_{L(\mathbf{s}, t)} \varphi(\mathbf{x}, t) \ dS \\
 &\quad \qquad + \frac{1}{d(\mathbf{s}, t)} \Big( \int_{L(\mathbf{s}, t)} \frac{d}{dt} \varphi(\mathbf{x}, t) \ dS + \frac{\partial_t d(\mathbf{s}, t)}{2} \varphi\rvert_{\mathbf{s}_-} + \frac{\partial_t d(\mathbf{s}, t)}{2} \varphi\rvert_{\mathbf{s}_+} \Big) \\
 &\quad = - \frac{\partial_t d(\mathbf{s}, t)}{d(\mathbf{s}, t)}  \varphi_\Gamma + \frac{1}{d(\mathbf{s}, t)} \int_{L(\mathbf{s}, t)} \frac{d}{dt} \varphi(\mathbf{x}, t) \ dS + \frac{\partial_t d(\mathbf{s}, t) }{d(\mathbf{s}, t)} \frac{\varphi\rvert_{\mathbf{s}_-} + \varphi\rvert_{\mathbf{s}_+}}{2}.
\end{align*}
Rearranging the terms results in (\ref{eq:timeintegral}).
\end{proof}
\end{lemma}

\begin{remark}
A similar result can be derived for the divergence of a vector-valued quantity $\bm{\psi} \in \left(C^1(\Omega_f^t)\right)^n$ that is split into tangential and normal part. It reads
\begin{equation}
\begin{split}
\label{eq:spaceintegral}
&\frac{1}{d(\mathbf{s}, t)} \int_{L(\mathbf{s}, t)} \diver \bm{ \psi }(\mathbf{x}, t) \ dS = \diver_\tau \bm{ \psi }_\Gamma(\mathbf{s}, t) \\ & \qquad
 + \frac{ \nabla_\tau d(\mathbf{s}, t) }{ d(\mathbf{s}, t) } \cdot \Big( \bm{\psi}_\Gamma(\mathbf{s}, t) - \mean{ \mathbf{T} \bm{\psi}(\mathbf{s}, t) } \Big)
  - \frac{1}{d(\mathbf{s}, t)} \jump{ \bm{\psi} \cdot \mathbf{n} }.
\end{split}
\end{equation}
\end{remark}

\begin{assumption}
\label{ass:integral}
Let us assume that there is no fluid exchange between bulk and fracture domain caused by the change of aperture, i.e., $\partial_t d(\mathbf{s}, t) (\phi S\rvert_{\mathbf{s}_\pm}) = 0$. Similarly, we will assume that $\mean{ \mathbf{T} \mathbf{v}^f(\mathbf{s}, t) }$ is small and neglect this term.
\end{assumption}

Then, for $\varphi = \phi^f S^f$ and $\bm{\psi} = \mathbf{v}^f$ (and similarly for $\bm{\psi} = \mathbf{F}^f(S^f, \mathbf{v}^f)$), using that $\phi^f$ is constant along $L(\mathbf{s}, t)$, (\ref{eq:timeintegral}) and (\ref{eq:spaceintegral}) read
\begin{align}
\label{eq:asstimeintegral}
\begin{split}
\frac{1}{d(\mathbf{s}, t)} \int_{L(\mathbf{s}, t)} \frac{d}{dt} &\left(\phi^f S^f(\mathbf{x}, t)\right) \ dS \\
 &= \frac{d}{dt} \left(\phi^f S_\Gamma(\mathbf{s}, t)\right) + \frac{ \partial_t d(\mathbf{s}, t) }{ d(\mathbf{s}, t) } \left( \phi^f S_\Gamma(\mathbf{s}, t)\right),
 \end{split}
 \\
\label{eq:assspaceintegral}
\int_{L(\mathbf{s}, t)} \diver \mathbf{v}^f(\mathbf{x}, t) \ dS
 &= \diver_\tau \left( d(\mathbf{s}, t) \mathbf{v}_\Gamma(\mathbf{s}, t) \right)
  - \jump{ \mathbf{v}^f \cdot \mathbf{n} }.
\end{align}

Now, integrating each equation in (\ref{eq:fffgequidim}) for $i = f$ along the line segments $L(\mathbf{s}, t)$, we obtain a reduced model for the fracture in terms of the unknowns $S_\Gamma$, $P_\Gamma$ and $\mathbf{v}_\Gamma$. Integrating the saturation conservation equation (\ref{eq:fffgequidim}) for $i = f$ and using (\ref{eq:asstimeintegral}), (\ref{eq:assspaceintegral}) we obtain
\begin{align}
q^\Gamma_w &= \frac{1}{d(\mathbf{s}, t)} \int_{L(\mathbf{s}, t)} ( \phi^f S^f )_t + \diver \mathbf{F}^f(S^f, \mathbf{v}^f) \ dS \\
 \label{eq:derivsat}
 &= ( \phi^f S_\Gamma )_t + 
    \frac{ \partial_t d(\mathbf{s}, t) }{ d(\mathbf{s}, t) } \phi^f S_\Gamma \\
    & \qquad \qquad + \frac{1}{d(\mathbf{s}, t)} \left( \diver_\tau \left( d(\mathbf{s}, t) \mathbf{F}^f( S_\Gamma, \mathbf{v}_\Gamma ) \right) - \llbracket \mathbf{F}(S, \mathbf{v}) \cdot \mathbf{n} \rrbracket \right) + O(\bar d)
 \label{eq:derivsatend}
\end{align}
with $q^\Gamma_w \coloneqq \frac{1}{d(\mathbf{s}, t)} \int_{L(\mathbf{s}, t)} q_w \ dS$.

We decompose $\mathbf{K}^f$ into $\mathbf{K}^f = K^f_n \mathbf{N} + \mathbf{K}^f_\tau \mathbf{T}$ and define $\mathbf{g}_\tau = \mathbf{T} \mathbf{g}$. Further, we multiply the equation (\ref{eq:fffgequidim}b) for $i = f$ by $\mathbf{T}$ and $\mathbf{N}$. By averaging over the line segments we deduce by (\ref{eq:relativeperdmdefA})-(\ref{eq:relativeperdmdefC}) the relation
\begin{align}
 0 &= \frac{1}{d(\mathbf{s}, t)} \int_{L(\mathbf{s}, t)} \mathbf{T} \mathbf{v}^f + \lambda^f(S^f) \mathbf{T} \mathbf{K}^f \big(\nabla P^f - G^f(S^f) \mathbf{g}\big)\ dS \\
 \label{eq:derivvel}
 &= \mathbf{v}_\Gamma + \frac{ \lambda^f(S_\Gamma) \mathbf{K}^f_\tau }{ d(\mathbf{s}, t) } \big( \nabla_\tau ( d(\mathbf{s}, t) P_\Gamma ) - G^f(S_\Gamma) \mathbf{g}_\tau \big) + O(\bar d)
\end{align}
and
\begin{align}
 0 &= \frac{1}{d(\mathbf{s}, t)} \int_{L(\mathbf{s}, t)} \mathbf{N} \mathbf{v}^f + \lambda^f(S^f) \mathbf{N} \mathbf{K}^f ( \nabla P^f - G^f(S^f) \mathbf{g}) \ dS \\
 \label{eq:derivcond}
 &= \mean{\mathbf{v} \cdot \mathbf{n}} + \lambda^f(S_\Gamma) K^f_n \left( - \frac{ \llbracket P \rrbracket }{d(\mathbf{s}, t)} - G^f(S_\Gamma) (\mathbf{g} \cdot \mathbf{n}) \right) + O(\bar d).
\end{align}
In (\ref{eq:derivcond}), we use that $\frac{1}{d(\mathbf{s}, t)} \int_{L(\mathbf{s}, t)} \mathbf{N} \mathbf{v}^f = \mean{\mathbf{v} \cdot \mathbf{n}} + O(\bar d^2)$. \\

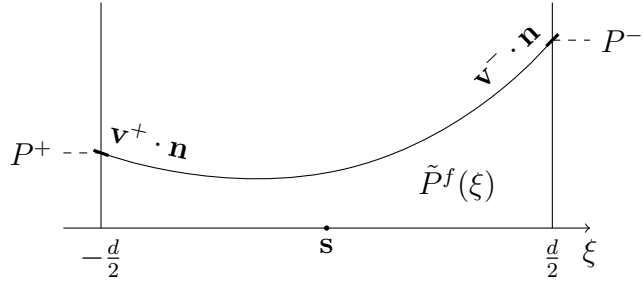
\begin{figure}[ht]
 \centering
  \begin{tikzpicture}
  \draw[thin] (-3,0) node [anchor=north]{$-\frac{d}{2}$} -- (-3,3);
  \draw[thin] (3,0) node [anchor=north]{$\frac{d}{2}$} -- (3,3);
  \draw plot [smooth, tension=1] coordinates {(-3,1) (0.3,0.8) (3,2.5)};
  \fill (0,0) circle (1pt) node[anchor=north] {$\mathbf{s}$};
  \node[anchor=south east] at (2.4,0.2) {$\tilde P^f(\xi)$};
  \draw [very thick] (-3.1,1.035) node [anchor=south west, rotate=-13]{$\mathbf{v}^+ \cdot \mathbf{n}$} -- (-2.9,0.965);
  \draw [dashed] (-3.5,1) node [anchor=east]{$P^+$} -- (-3,1);
  \draw [very thick] (3.08,2.58) node [anchor=south east, rotate=44]{$\mathbf{v}^- \cdot \mathbf{n}$} -- (2.92,2.42);
  \draw [dashed] (3,2.5) -- (3.5,2.5) node [anchor=west]{$P^-$};
  \draw[->] (-3.5,0) -- (3.5,0) node[anchor=north] {$\xi$};
 \end{tikzpicture}
 \caption{Within the fractures the pressure along the line segments is approximated by a second-order polynomial.}
 \label{fig:pressurePolynomial}
\end{figure}

We approximate $P^f$ by a quadratic polynomial $\tilde P^f(\xi) = a \xi^2 + b \xi + c,\ \xi \in \left[-\frac{d}{2}, \frac{d}{2}\right],$ along the line segment $L(\mathbf{s}, t)$. This polynomial is supposed to satisfying the four continuity conditions (\ref{eq:couplcondsB}), (\ref{eq:couplcondsC}) at the boundaries (compare Fig.\ \ref{fig:pressurePolynomial}) where one can be eliminated using (\ref{eq:derivcond}). Averaging as in (\ref{eq:avgquansC}), we obtain the relation
\begin{align}
    \label{eq:pfaverage}
    P_\Gamma &= \frac{1}{d(\mathbf{s}, t)} \int_{-\frac{d}{2}}^{\frac{d}{2}} \tilde P^f(\xi)\ d\xi  + O(\bar d) \\
    &= \mean{P} - \frac{d} {12 \lambda^f(S_\Gamma) \mathbf{K}_n^f} \jump{\mathbf{v} \cdot \mathbf{n}} + O(\bar d).
\end{align}

Combining equation (\ref{eq:derivcond}) and (\ref{eq:pfaverage}) leads for $\alpha \in \{+,-\}$ to
\begin{align}
 \begin{split}
 \label{eq:cond2side}
 \mathbf{v}^\alpha \cdot \mathbf{n}^\alpha = - \lambda^f(S_\Gamma) K^f_n \Bigg( \frac{P_\Gamma - P^\alpha}{d / 2} + \frac{ P_\Gamma - \mean{P} }{d / 4} \\
 -\ G^f(S_\Gamma) (\mathbf{g} \cdot \mathbf{n}^\alpha) \Bigg) + O(\bar d).
 \end{split}
\end{align}
This identity provides a coupling condition between the bulk and the fracture problem, cf. (\ref{eq:ffflowdimcoupling}) below.

From the divergence constraint in (\ref{eq:fffgequidim}) we obtain using (\ref{eq:spaceintegral}) the identity
\begin{align}
  \int_{L(\mathbf{s}, t)} \diver( \mathbf{v}^f ) \ dS &= \diver_\tau(d(\mathbf{s}, t) \mathbf{v}_\Gamma) - \llbracket \mathbf{v} \cdot \mathbf{n} \rrbracket \\
  &= \int_{L(\mathbf{s}, t)} (q_w + q_{nw})\ dS = d(\mathbf{s}, t) \left( q_w^\Gamma + q_{nw}^\Gamma \right),
\end{align}
with $q_{nw}^\Gamma \coloneqq \frac{1}{d(\mathbf{s}, t)} \int_{L(\mathbf{s},t)} q_{nw} \ dS$.

The coupling conditions for the mass balance (\ref{eq:couplcondsA}) can be rewritten as
\begin{align}
 \label{eq:derivcoupsat}
 \begin{split}
 &\left(\mathbf{F}(S^\alpha(\mathbf{s}_\alpha, t), \mathbf{v}^\alpha(\mathbf{s}_\alpha, t)) \cdot \mathbf{n}^\alpha \right) \rvert_{\gamma_\alpha} \\
 &= \left(\mathbf{F}^f(S^f(\mathbf{s}_\alpha, t), \mathbf{v}^f(\mathbf{s}_\alpha, t)) \cdot \mathbf{n}^\alpha\right) \rvert_{\gamma_\alpha} \\
 &= \left(\mathbf{F}^f(S_\Gamma(\mathbf{s}, t), \mathbf{v}^\alpha(\mathbf{s}_\alpha, t)) \cdot \mathbf{n}^\alpha \right)\rvert_{\gamma_\alpha} + O(\bar d)
 \end{split}
\end{align}
for $\alpha \in \{+,-\}$.

\subsection{The Reduced Model}
\label{sec:reducedModel}

Next, we summarize the derived equations (\ref{eq:derivsat})-(\ref{eq:derivsatend}), (\ref{eq:derivvel}), (\ref{eq:cond2side}) and (\ref{eq:derivcoupsat}) neglecting the $O(\bar d)$-terms.
For each $t  \in (0,T)$ let as before $\Gamma(t) \subset \Omega$ be a given family of hypersurfaces and $d(\mathbf{s}, t)$ the corresponding aperture in $\mathbf{s} \in \Gamma(t)$ and $D(t) \coloneqq \Omega \setminus \bar \Gamma(t)$.

The wetting fluid saturation $S: D_t \to [0,1]$, the global pressure $P: D_t \to \mathbb{R}$ and the total velocity $\mathbf{v}: D_t \to \mathbb{R}^n$ in the bulk medium satisfy

\begin{align}
 \label{eq:ffflowdim}
 \left.
 \begin{aligned}
  (\phi S)_t + \diver \mathbf{F}(S, \mathbf{v}) &= q_w, \\
  \mathbf{v} + \lambda(S) \mathbf{K} ( \nabla P - G(S) \mathbf{g} ) &= \mathbf{0}, \\
  \diver(\mathbf{v}) &= q_w + q_{nw}
 \end{aligned}
 \qquad \right. \quad
 \begin{gathered}
  \text{in } D_t,
 \end{gathered}
\end{align}
where $\mathbf{F}(S, \mathbf{v})$ and $G(S)$ are defined as for (\ref{eq:fffg}).

On $\Gamma_t$ we search for the reduced quantities $S_\Gamma: \Gamma_t \to [0,1],\ P_\Gamma: \Gamma_t \to \mathbb{R}$ and $\mathbf{v}_\Gamma: \Gamma_t \to \mathbb{R}^n$ satisfying

\begin{align}
 \label{eq:ffffraclowdim}
 \begin{aligned}
  (d \phi^f S_\Gamma)_t + \diver_\tau \left( d \mathbf{F}^f(S_\Gamma, \mathbf{v}_\Gamma) \right) &= \llbracket \mathbf{F}(S, \mathbf{v}) \cdot \mathbf{n} \rrbracket + d q_w^{\Gamma}, \\ 
  d\mathbf{v}_\Gamma + \lambda^f(S_\Gamma) \mathbf{K}^f_\tau \left( \nabla_\tau (d P_\Gamma) - G^f(S_\Gamma) \mathbf{g}_\tau \right) &= \mathbf{0}, \qquad \qquad \qquad \qquad \text{in } \Gamma_t. \\
  \diver_\tau(d \mathbf{v}_\Gamma) &= \llbracket \mathbf{v} \cdot \mathbf{n} \rrbracket + d ( q_w^{\Gamma} + q_{nw}^{\Gamma} )
 \end{aligned}
\end{align}
The systems (\ref{eq:ffflowdim}), (\ref{eq:ffffraclowdim}) are closed at the hypersurface $\Gamma(t)$ by
\begin{align}
 \begin{aligned}
 \label{eq:ffflowdimcoupling}
 &\mathbf{F}(S^\alpha, \mathbf{v}^\alpha) \cdot \mathbf{n}^\alpha = \mathbf{F}^f(S_\Gamma, \mathbf{v}^\alpha) \cdot \mathbf{n}^\alpha, \\
 &\mathbf{v}^\alpha \cdot \mathbf{n}^\alpha = - \lambda^f(S_\Gamma) K^f_n \left( \frac{P_\Gamma - P^\alpha}{d/2} + \frac{ P_\Gamma - \mean{P} }{d / 4} - G^f(S_\Gamma) (\mathbf{g} \cdot \mathbf{n}^\alpha) \right)
 \end{aligned}
\end{align}
on $\Gamma(t) \times \{t\}, t \in (0,T), \alpha \in \{+,-\}$.

\begin{remark}
(i) In case of $S \equiv S_\Gamma \equiv 0$ and $\mathbf{g} = \mathbf{0}$ the model (\ref{eq:ffflowdim})-(\ref{eq:ffflowdimcoupling}) reduces to the single-phase case as described in \cite{VincentMartin2004}.

(ii) For $d = const.$ we reconstruct the model as in \cite{Fumagalli2013}, but obtain a slightly different coupling condition. This is because during model deduction we replaced $G^f(S)\rvert_{\gamma_\pm}$ by $G^f(S_\Gamma)$ instead of $G^f(S^\pm)$. This choice appears consistent to us as it is used for all other terms depending on $S$.
\end{remark}

The equations (\ref{eq:ffflowdim}), (\ref{eq:ffffraclowdim}) describe a mixed-dimensional problem of hyper\-bolic-elliptic type in time-dependent domains. Special numerical methods and tools are required to solve such kind of problem.

In the next section we will describe our solution approach.
We follow a finite-volume approach because of the hyperbolic character of the saturation equation. As a novel contribution this approach is coupled with a moving-mesh concept that keeps track of the moving lower-dimensional domain.

\section{The Numerical Scheme}

A large variety of numerical methods have been proposed for mixed-dim\-ensional models (see e.g. \cite{Fumagalli2013, Berre2021, Antonietti2019}).

They can be classified in two categories: conforming and non-conforming methods. The non-conforming methods make use of independent discretisations of full dimensional and lower-dimensional sets, whereas conforming methods assume some kind of conformity of the lower-dimensional set to the bulk mesh. For instance, mesh elements of the lower-dimensional mesh coincide with facets of the bulk mesh. Non-conforming methods seem to suggest themselves for moving interfaces being much more easy to handle. However, the mutual geometrical relations are not trivial to sustain. If instead one is able to manage a moving mesh, the geo\-metrical relations become trivial. Therefore, we propose a moving-mesh method that permits the efficient tracking of lower-dimensional mesh facets and utilize the conformity in a finite-volume method that is just slightly enhanced by adopted fluxes at the lower-dimensional interface.

In the following a finite-volume discretization will be used for both the bulk and the fracture problem. The coupling between the two problems is incorporated by adopted fluxes at the edges of the codimension-1 interface which ensures in particular conservation of mass. A finite-volume moving-mesh (FVMM) scheme is implemented to keep track of the moving interface.

We will use a two-point flux approximation (TPFA) for the discretization of the pressure gradient in (\ref{eq:ffflowdim}), (\ref{eq:ffffraclowdim}). We use this simple approach because large effort is necessary to update cell stencils during re-meshing. It is known that a standard TPFA is not consistent on triangular grids and for anisotropic permeability tensors \cite{Eymard2000}. Therefore, we will propose an adaptation of the TPFA method that is consistent at least for isotropic permeabilities in Section \ref{subsec:tpfacirc} locating the pressure values at the cirumcenters.

Before we start with the description of the scheme, let us introduce some notation.

\subsection{Notation}

Let $0 = t_0 < t_1 < \dots < t_N = T, N \in \mathbb{N}$, be a series of time steps and $\mathcal{T} = \mathcal{T}(t)$ a conforming, time-dependent triangulation of $\Omega$ such that a subset of the facets of $\mathcal{T}$ coincides with a $(n-1)$-dimensional triangulation $\mathcal{T}_\Gamma = \mathcal{T}_\Gamma(t)$ of $\Gamma(t)$ for all $t \in [0,T]$.

\begin{wrapfigure}{r}{0.4\textwidth}
 \centering
  \begin{tikzpicture}[scale=2]

  \coordinate (A) at (-0.1,-0.1);
  \coordinate (B) at (0.7,0.6);
  \coordinate (C) at (1.5,1.5);
  \coordinate (D) at (0.0,1.2);
  \coordinate (E) at (0.5,1.7);
  \coordinate (F) at (0.8,-0.3);
  \coordinate (G) at (1.6,0.3);
  \draw (A) -- (D) -- (B);
  \draw (D) -- (E) -- (B);
  \draw (B) -- (E) -- (C);
  \draw (A) -- (F) -- (B);
  \draw (C) -- (G) -- (B);
  \draw (F) -- (G);
  \draw[Blue, very thick] (A) -- (B) -- (C);

  \node at (1.6,0) {$\mathcal{T}$};
  \draw[very thin] (1.15,0.45) -- (1.1,0.3) node[anchor=north]{$\mathcal{F}^I$};
  \draw[Blue, very thin] (B) -- (0.83,1.0) node[anchor=south]{$\mathcal{F}_\Gamma$};
  \fill[Blue] (A) circle (1pt);
  \fill[Blue] (B) circle (1pt);
  \fill[Blue] (C) circle (1pt) node[anchor=south, yshift=0.1cm] {$\mathcal{T}_\Gamma (= \mathcal{F}^\Gamma)$};
  
  \draw[thick, ->] (0.3,0.25) -- +(-0.18, 0.21) node [anchor=north,yshift=-0.05cm] {\tiny$\mathbf{n}$};
  \node at (0.25,0.65) {$K^-$};
  \node at (0.52,0.07) {$K^+$};
 \end{tikzpicture}
 \caption{Visualization of the geometrical notation for $n=2$.}
 \label{fig:grid}
\end{wrapfigure}
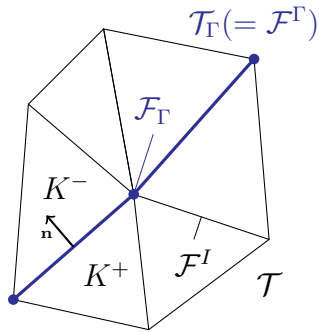

Denote the set of facets of the triangulation $\mathcal{T}$ by $\mathcal{F}$ and the facets of $\mathcal{T}_\Gamma$ by $\mathcal{F}_\Gamma$, respectively. We distinguish between inner facets $\mathcal{F}^I$ and those of $\mathcal{F}$ which coincide with the lower-dimensional mesh $\mathcal{F}^\Gamma$. Further, let us denote by $\mathcal{F}(K)$ the facets of a cell $K$ and denote by $\mathbf{n}$ the outer normal to this facet multiplied by the facet's area. The outer normal defines an inner, positive side of the facet and an outer, negative side. The adjacent cells we denote therefore by $K^+$ and $K^-$. We assume that the vertices of $\mathcal{T}$ move linearly in time in each time interval $(t_n, t_{n+1})$. That is, all $K = K(t) \in \mathcal{T}(t)$ are given by $K(t) \coloneqq \text{conv}( \mathbf{p}_0(t), \dots, \mathbf{p}_n(t) )$, where $\mathbf{p}_i(t) \in \mathbb{R}^n$ moves with speed $\mathbf{s}_i \in \mathbb{R}^n$ according to
\begin{align}
 \mathbf{p}_i(t) = \mathbf{p}_i(0) + t \mathbf{s}_i, \qquad i=0, \dots, n.
\end{align}

\subsection{The Finite-Volume Moving-Mesh Method}

In this section, we describe the FVMM method that is used for handling the propagation of fractures.
Usually, moving-mesh methods are used to minimize artificial diffusion solving, e.g., hyperbolic problems tracking discontinuities of the solution \cite{Harten83}. In contrast, we use it to track a lower-dimensional interface.

The FVMM methods require an additional geometrical flux within the finite-volume formulation. This geometrical flux accounts for the mass flux across moving edges \cite{Chalons2017}. Because it is only an additional flux in the finite-volume update step, it works without any re-meshing and projection as long as the triangulation is not adapted.
With this method one is able to move the lower-dimensional fracture facets without moving mass in the solution of the surrounding mesh.

Let us consider a single interval $(t_n, t_{n+1})$. We assume that $\mathcal{T}(t)$ does not degenerate for $t \in (t_n, t_{n+1})$, i.e., there is a $c > 0$ such that for all $t \in (t_n, t_{n+1})$ $|K| > c\ \forall K \in \mathcal{T}(t)$. If this condition is not satisfied, the triangulation has to be adapted first, as described in Section \ref{subsec:implementation}.

Now, consider the space-time cell
\begin{align}
 K_{st} = \{ (\mathbf{x},t)\ |\ \mathbf{x} \in K(t),\ t_n \leq t \leq t_{n+1} \}.
\end{align}

Integrating (\ref{eq:ffflowdim}a) over $K_{st}$ and using Reynolds' transport theorem we compute
\begin{align}
\int_{K_{st}} (\phi S)_t & + \diver \mathbf{F}(S, \mathbf{v}) \ d(\mathbf{x},t) \\
\label{eq:reynolds1}
=& \int_{K(t_{n+1})} \phi S(\cdot, t_{n+1}) \ d\mathbf{x} - \int_{K(t_n)} \phi S(\cdot, t_n) \ d\mathbf{x}  \\
\label{eq:reynolds2}
&+ \int_{t_n}^{t_{n+1}} \int_{\partial K(t)} \big(\mathbf{F}(S, \mathbf{v}) - (\phi S) \mathbf{s} \big) \cdot \mathbf{n} \ dS\ dt \\
&= \int_{t_n}^{t_{n+1}} \int_{K(t)} q_w \ dS\ dt.
\end{align}

Here, $\mathbf{s}: K(t) \to \mathbb{R}^n$ is the speed of a point
\begin{align}
 \mathbf{x}(t) = \mathbf{p}_0(t) + \sum_{i=1}^n \lambda_i \big( \mathbf{p}_i(t) - \mathbf{p}_0(t) \big),
\end{align}
with $\lambda_i \in [0,1],\ i=1, \dots, n$, i.e.,
\begin{align}
 \mathbf{s}(\mathbf{x}) = \mathbf{s}_0 + \sum_{i=1}^n \lambda_i (\mathbf{s}_i - \mathbf{s}_0).
\end{align}
Further, $\mathbf{n}: \partial K(t) \to \mathcal{S}^{n-1}$ denotes the unit outer normal at the boundary of $K(t)$.
Let us define the finite-volume ansatz-space of cell-wise constant functions by
\begin{align}
 \mathcal{S}_{h}(\mathcal{T}) \coloneqq \{ v \in L^2(\mathbb{R}^d)\ |\ v_K \coloneqq v\rvert_K \in P_0(K)\ \forall K \in \mathcal{T} \}.
\end{align}
We choose two discrete representatives of the saturation $S^n, S^{n+1} \in \mathcal{S}_{h}(\mathcal{T})$ and define $S^0$ by $S^0_K \coloneqq \frac{1}{|K|} \int_K S(\cdot, 0) \ d\mathbf{x}$.
Now, using an implicit Euler time-stepping for $(\ref{eq:reynolds1}),(\ref{eq:reynolds2})$ and dividing by $| K(t_{n+1}) |$ we obtain the successive definition of $S^n \to S^{n+1}$ given by
\begin{align}
 \phi S_K^{n+1} & - \phi S_K^n \frac{ | K(t_n) | }{ | K(t_{n+1}) | } \\
 &+ (t_{n+1} - t_n) \sum_{\mathcal{F}\big(K(t_{n+1})\big)} \Big[ g(S^{n+1}_{K^+}, S^{n+1}_{K^-}, \mathbf{v}) + h(S^{n+1}_{K^+}, S^{n+1}_{K^-}) \Big] \\
 &= \frac{t_{n+1} - t_n}{| K(t_{n+1}) |} \int_{K(t_{n+1})} q_w \ d\mathbf{x}.
\end{align}

Here, $g(\cdot, \cdot, \mathbf{v})$ is a numerical flux that is consistent with the flux function $\mathbf{F}(S, \mathbf{v}) \cdot \mathbf{n}$, whereas the numerical flux $h(\cdot, \cdot)$ has to be consistent with the flux function $-(\phi S) \mathbf{s} \cdot \mathbf{n}$. Suitable choices are, for instance, a Lax-Friedrichs or Godunov-type flux for $g$ and an upwind flux for $h$.
Because the flux function $\mathbf{F}(S, \mathbf{v})$ in our model is non-monotone in the first argument, we will use the Godunov flux that results from an exact solution of the Riemann problem, c.f., Section \ref{subsec:godunov}.
The fluxes at the boundary have to be adopted according to the boundary conditions.

With the additional geometrical flux and volume term the conservative quantity is not transported in space although vertices of the triangulation move. We show in Section \ref{subsec:case1} numerical examples how the scheme performs for constant initial data.

Of course, this scheme can only be applied as long as no cell degenerates. In order to prevent degeneration and improve mesh quality we use the re-meshing techniques provided by \dunemmesh\ \cite{Mmesh}, see Section \ref{subsec:implementation} for details.

We can use this FVMM method for both the bulk and the lower-dim\-ensional domain. To couple the schemes, the fluxes at the inner fracture facets $\mathcal{F}^\Gamma$ have to incorporate condition (\ref{eq:ffflowdimcoupling}).

\subsection{The Generalized Godunov flux}
\label{subsec:godunov}
We shall describe our choice for the numerical flux $g(\cdot, \cdot, \mathbf{v})$ in more detail. The flux has to be generalized for a discontinuous flux function as the physical properties (in particular $\mathbf{K}$) might vary in space.

We start with the formulation of the Godunov flux \cite{kroener}
\begin{align}
    g(S^+, S^-, \mathbf{v}) = \begin{cases}
     \min\limits_{S^+ < S < S^-} \mathbf{F}(S, \mathbf{v}) \cdot \mathbf{n},& \text{if } S^+ \leq S^-, \\
     \max\limits_{S^- < S < S^+} \mathbf{F}(S, \mathbf{v}) \cdot \mathbf{n},& \text{if } S^- < S^+.
    \end{cases}
\end{align}
Generalizing this flux for a discontinuous flux function we use the identity
\begin{align}
    g^*(S^+, S^-, \mathbf{v}) \coloneqq g^+(S^+, S^*, \mathbf{v}) = g^-(S^*, S^-, \mathbf{v}) \text{ for some } S^* \in [0,1].
\end{align}
This generalization has to be applied on all facets where the physical parameters of the adjacent cells are distinct, in particular between bulk and fracture domain.

For a quadratic material law $k_w(S) = S^2$ and $k_{nw}(S) = (1-S)^2$ we can compute explicit formulas, that can be implemented efficiently.
Therefore, we use the fact that $\mathbf{F}(S, \mathbf{v})$ has a single extremum (compare Fig. \ref{fig:fracflow}). This fact can also be exploited to deduce an explicit formula for the generalized flux.

\begin{figure}[ht]
\begin{minipage}{0.49\textwidth}
 \centering
  \begin{tikzpicture}[scale=0.7]
  \begin{axis}[
    axis lines = left,
    xlabel = {$S$},
    ylabel = {$k_\alpha(S)$},
  ]
  
  \addplot [
    domain=0:1, 
    samples=101, 
    color=black,
    style=thick
  ]
  {x^2};
  
  \addplot [
    domain=0:1, 
    samples=101, 
    color=gray,
    style=thick
  ]
  {x};

  \addplot [
    domain=0:1, 
    samples=101, 
    color=black,
    style=thick, dashed
    ]
    {(1-x)^2};
    
  \addplot [
    domain=0:1, 
    samples=101, 
    color=gray,
    style=thick, dashed
    ]
    {(1-x)};
   
  \end{axis}
 \end{tikzpicture}
 \caption{Relative permeability functions for two (linear/quadratic) choices for $k_\alpha$. The continuous curve shows $k_w$, the dashed one shows $k_{nw}$.}
 \label{fig:relperm}
\end{minipage}
\hfill
\begin{minipage}{0.49\textwidth}
 \centering
  \begin{tikzpicture}[scale=0.7]
  \begin{axis}[
    axis lines = left,
    xlabel = {$S$},
    ylabel = {$\mathbf{F}(S,\mathbf{v}) \cdot \mathbf{n}$},
    legend pos=north west,
    legend cell align={left}
  ]

  \addplot [
    domain=0:1, 
    samples=101, 
    color=gray,
    style=thick
  ]
  {x / ( x + (1-x) ) * ( 2 - (1-x)*10 )};
  
  \addlegendentry{linear}
  
  \addplot [
    domain=0:1, 
    samples=101, 
    color=black,
    style=thick
  ]
  {x^2 / ( x^2 + (1-x)^2 ) * ( 2 - (1-x)^2*10};
  
  \addlegendentry{quadratic}
  
  \end{axis}
 \end{tikzpicture}
 \caption{The non-monotone shape of the flux function $\mathbf{F}(\cdot, \mathbf{v}) \cdot \mathbf{n}$. Here, $k_\alpha$ is linear/quadratic, $\mathbf{v} \cdot \mathbf{n}=2$ and $\mathbf{n} \mathbf{K} (\rho_n - \rho_w) \mathbf{g} = 10$ .}
 \label{fig:fracflow}
\end{minipage}
\end{figure}

\subsection{Circum-centered Two-Point Flux Approximation}
\label{subsec:tpfacirc}
We continue with the discretization of the elliptic part (\ref{eq:ffflowdim}b)-(\ref{eq:ffflowdim}c) governing bulk pressure and velocity. The derivation shows the method for the bulk problem in $D(t)$, but it is similar for the problem on the interface $\Gamma(t)$.

The finite-volume approach for the divergence constraint (\ref{eq:fffg}c) reads
\[
 \sum_{F \in \mathcal{F}\left(K(t)\right)} v_F \ dS = \int_{K(t)} ( q_w + q_{nw} ) \ d\mathbf{x},
\]
where $v_F$ is a suitable approximation of $\mathbf{v} \cdot \mathbf{n}$ on $F$ and $K(t) \in \mathcal{T}(t)$. A simple choice for $v_F$ is the two-point flux approximation derived from equation (\ref{eq:fffg}b), see \cite{Eymard2000}.
Including the gravity term it reads
\begin{align*}
 \label{eq:tpfa}
 v_F &\coloneqq - \mathbb{T}_F \big( P_{K^-} - P_{K^+} - \mathbb{G}_F \big),
\end{align*}
where the transmissibility $\mathbb{T}_F$ is defined by
\begin{align*}
 \mathbb{T}_F \coloneqq \frac{\mathbb{T}_{K^+} \mathbb{T}_{K^-}}{\mathbb{T}_{K^+} + \mathbb{T}_{K^-}} \qquad \text{with} \qquad
 \mathbb{T}_i \coloneqq \lambda(S_i) \frac{\mathbf{d}_i \mathbf{K}_i \mathbf{d}_i}{\| \mathbf{d}_i \|_2^3}, \ i \in \{ K^+, K^- \}.
\end{align*}
Here, $\mathbf{d}_i \coloneqq \mathbf{m}_{F} - \mathbf{m}_i,\ i \in \{ K^+, K^- \}$, is the distance vector between the center of the facet $F$ and the cell centers. $\mathbf{K}_i$ for $i \in \{ K^+, K^- \}$ denotes the restriction of $\mathbf{K}$ to $K^+$ and $K^-$. For consistency of the scheme it is necessary that $\mathbf{m}_{K^+}, \mathbf{m}_{K^-}, \mathbf{m}_F$ are the cirumcenters of $K^+, K^-, F$.
The gravitational influence $\mathbb{G}_F$ is defined by
\begin{align*}
 \mathbb{G}_F \coloneqq \mathbb{G}_{K^+} - \mathbb{G}_{K^-}, \qquad
 \mathbb{G}_i \coloneqq G(S_i)\ (\mathbf{d}_i \cdot \mathbf{g}), \ i \in \{ K^+, K^- \}.
\end{align*}

The choices for the transmissibilities have been made such that
\begin{align}
    - \mathbb{T}_{K^+} \big( P^* - P_{K^+} - \mathbb{G}_{K^+} \big) = \mathbb{T}_{K^-} \big( P^* - P_{K^-} - \mathbb{G}_{K^-} \big)
\end{align}
holds for some intermediate pressure value $P^*$ that can be eliminated.

\newcommand{\tplus}{K^+}
\newcommand{\tminus}{K^-}

At facets that coincide with a lower-dimensional fracture element $K_\Gamma \in \mathcal{T}_\Gamma$ we include the coupling conditions (\ref{eq:ffflowdimcoupling}) of the reduced model. Therefore, we introduce intermediate pressure values $P \rvert_{\gamma_+}$ and $P \rvert_{\gamma_-}$ at the boundaries of the bulk medium next to the fracture.
Conditions for the intermediate pressure values can be stated by
\begin{align}
 v_F \rvert_{\gamma_+} &= - \mathbb{T}_{\tplus} \big( P \rvert_{\gamma_+} - P_{\tplus} - \mathbb{G}_{\tplus} \big), \\
 v_F \rvert_{\gamma_-} &= - \mathbb{T}_{\tminus} \big( P \rvert_{\gamma_-} - P_{\tminus} - \mathbb{G}_{\tminus} \big).
\end{align}
Then, the coupling conditions in (\ref{eq:ffflowdimcoupling}) can be used to eliminate the intermediate values. Defining
\begin{align}
    \mathbb{T}_\Gamma \coloneqq \frac{2}{d} \lambda^f(S_{T_\Gamma})K^f_n, \qquad
    \mathbb{G}_\Gamma \coloneqq - \frac{d}{2} G^f(S_{K_\Gamma})\ (\mathbf{n} \cdot \mathbf{g})
\end{align}
we obtain
\begin{align}
 \label{eg:couplvel}
 v_F \rvert_{\gamma_+} &= R
 \begin{pmatrix}
    3\mathbb{T}_\Gamma + 2 \mathbb{T}_{\tminus} \\ 
    - 3\mathbb{T}_\Gamma - 3\mathbb{T}_{\tminus} \\
    \mathbb{T}_{\tminus}
 \end{pmatrix}
 \cdot
 \begin{pmatrix}
    P_{\tplus} + \mathbb{G}_{\tplus} \\ 
    P_\Gamma + \mathbb{G}_\Gamma \\
    P_{\tminus} + \mathbb{G}_{\tminus} + 2 \mathbb{G}_\Gamma
 \end{pmatrix}
\end{align}
where
\begin{align}
R \coloneqq \frac{ \mathbb{T}_{\tplus} \mathbb{T}_\Gamma}{\mathbb{T}_{\tplus} \mathbb{T}_{\tminus} + 3 {\mathbb{T}_\Gamma}^2 + 2\mathbb{T}_\Gamma (\mathbb{T}_{\tplus} + \mathbb{T}_{\tminus} )}.
\end{align}

\begin{remark}
The presented discretization is consistent for isotropic intrinsic permeabilities if we use circumcenters for $\mathbf{m}_{K^+}, \mathbf{m}_{K^-}, \mathbf{m}_F$, and therefore locate the pressure values at the circumcenters of the tetrahedral cells \cite{Eymard2000}.
This is still valid for the coupling to the fracture network as the circumcenters of the lower-dimensional mesh elements are located at the orthogonal connection line of the circumcenters of the two adjacent bulk cells.
\end{remark}

We have shown that this scheme produces quite comparative results in a recent benchmark study \cite{Berre2021}.

\clearpage
\subsection{The Complete FVMM Algorithm}
Let us summarize the complete scheme for bulk and fracture domain. It reads as follows.

For each time step $t_n$ with given $S^{n} \in \mathcal{S}_{h}(\mathcal{T})$ and $(S_\Gamma)^{n} \in \mathcal{S}_{h}(\mathcal{T}_\Gamma)$, we solve the following system implicitly for $(S^{n+1}, P^{n+1}) \in \mathcal{S}_{h}(\mathcal{T}) \times \mathcal{S}_{h}(\mathcal{T})$, $\left((S_\Gamma)^{n+1}, (P_\Gamma)^{n+1}\right) \in \mathcal{S}_{h}(\mathcal{T}_\Gamma) \times \mathcal{S}_{h}(\mathcal{T}_\Gamma)$ defined by

\begin{align}
 \begin{split}
 & \frac{ \phi S_K^{n+1} - \phi S_K^n \frac{ | K(t_n) | }{ | K(t_{n+1}) | } } { t_{n+1} - t_n }
  + \sum_{F \in \mathcal{F}^I} \big( g(S_{K^+}^{n+1}, S_{K^-}^{n+1}, v_F^{n+1}\mathbf{n} ) + h(S_{K^+}^{n+1}, S_{K^-}^{n+1}) \big) \\
 & \qquad \qquad
  + \sum_{F \in \mathcal{F}^\Gamma} g(S_{K^+}^{n+1}, S_{K_\Gamma}^{n+1}, v_F\rvert_{\gamma_+}^{n+1} \mathbf{n})
 = \int_{ K(t_{n+1}) } q_w(\mathbf{x}, t_{n+1}) \ d\mathbf{x}
 \end{split} \\
 & \sum_{F \in \mathcal{F}^I} v_F^{n+1}\mathbf{n} 
 + \sum_{F \in \mathcal{F}^\Gamma} v_F\rvert_{\gamma_+}^{n+1}\mathbf{n}
 = \int_{ K(t_{n+1}) } \left(q_w(\mathbf{x}, t_{n+1}) + q_{nw}(\mathbf{x}, t_{n+1})\right) \ d\mathbf{x}
\end{align}
for all $K(t) \in \mathcal{T}(t)$ and
\begin{align}
 \begin{split}
 & \frac{ d^{n+1} \phi^\Gamma S_{K_\Gamma}^{n+1} - d^{n} \phi^\Gamma S_{K_\Gamma}^n \frac{ | {K_\Gamma}(t_n) | }{ | {K_\Gamma}(t_{n+1}) | } } { t_{n+1} - t_n } \\[2mm]
 & \qquad \qquad 
  + \sum_{F \in \mathcal{F}^I_\Gamma} \big( d^{n+1}  g(S_{K^+_\Gamma}^{n+1}, S_{K^-_\Gamma}^{n+1}, \mathbf{v}_{F_\Gamma}\mathbf{n}_\Gamma)
  + h(S_{K^+_\Gamma}^{n+1}, S_{K^-_\Gamma}^{n+1}) \big) \\[2mm]
  & \qquad =
 g(S_{K_\Gamma}^{n+1}, S_{K^+_\Gamma}^{n+1}, v_F\rvert_{\gamma_+}^{n+1}\mathbf{n} ) + g(S_{K_\Gamma}^{n+1}, S_{K^-_\Gamma}^{n+1}, v_F\rvert_{\gamma_-}^{n+1}\mathbf{n}) \\[3mm]
 & \qquad \qquad
 + d^{n+1} \int_{ K_\Gamma(t_{n+1}) } q_w^\Gamma(\mathbf{x}, t_{n+1}) \ d\mathbf{x},
 \end{split} \\[2mm]
 \begin{split}
 & \sum_{F \in \mathcal{F}^I_\Gamma} d^{n+1} \mathbf{v}_{F_\Gamma}
 = v_F\rvert_{\gamma_+}^{n+1}\mathbf{n} + v_F\rvert_{\gamma_-}^{n+1}\mathbf{n} \\
 & \qquad \qquad + d^{n+1} \int_{ K_\Gamma(t_{n+1}) } \left(q_w^\Gamma(\mathbf{x}, t_{n+1}) + q_{nw}^\Gamma(\mathbf{x}, t_{n+1})\right) \ d\mathbf{x}
 \end{split}
\end{align}
for all $K_\Gamma(t) \in \mathcal{T}_\Gamma(t)$.

Here, $d^{n} = d(\cdot, t^{n})$. For the purpose of readability we neglect the boundary terms in the formulation above. At Dirichlet boundaries, the outer values $S_{K^-}^{n+1}$ and $P_{K^-}^{n+1}$ (respectively $S_{K^-_\Gamma}^{n+1}$ and $P_{K^-_\Gamma}^{n+1}$) have to be replaced by the Dirichlet boundary value. At Neumann boundaries, the corresponding normal fluxes $g$ and $v_F$ can be replaced directly by the Neumann boundary flux.

\subsection{Implementation}
\label{subsec:implementation}

We implemented our method within the software framework DUNE \cite{dune} on the basis of the discretization module \dumux \cite{dumux} and the grid implementation \dunemmesh\ \cite{Mmesh}. The grid implementation \dunemmesh\ is a new development and is essential for both the mixed-dimensional discretization and the moving-mesh method. It is a grid wrapper of CGAL \cite{cgal} triangulations in 2D and 3D and can export a pre-described set of facets as a separate network grid. One of the main advantages of the strong coupling of the two grids is the simultaneous re-meshing of bulk and interface grid for dimension $n=2$. An open-source release of \dunemmesh\ is available \cite{Mmesh}.
The re-meshing feature of \dunemmesh\ is able to insert and remove arbitrary vertices at any time and assists in projecting unknowns. 

\begin{wrapfigure}{r}{0.45\textwidth}
    \centering
	\includegraphics[width=0.18\textwidth]{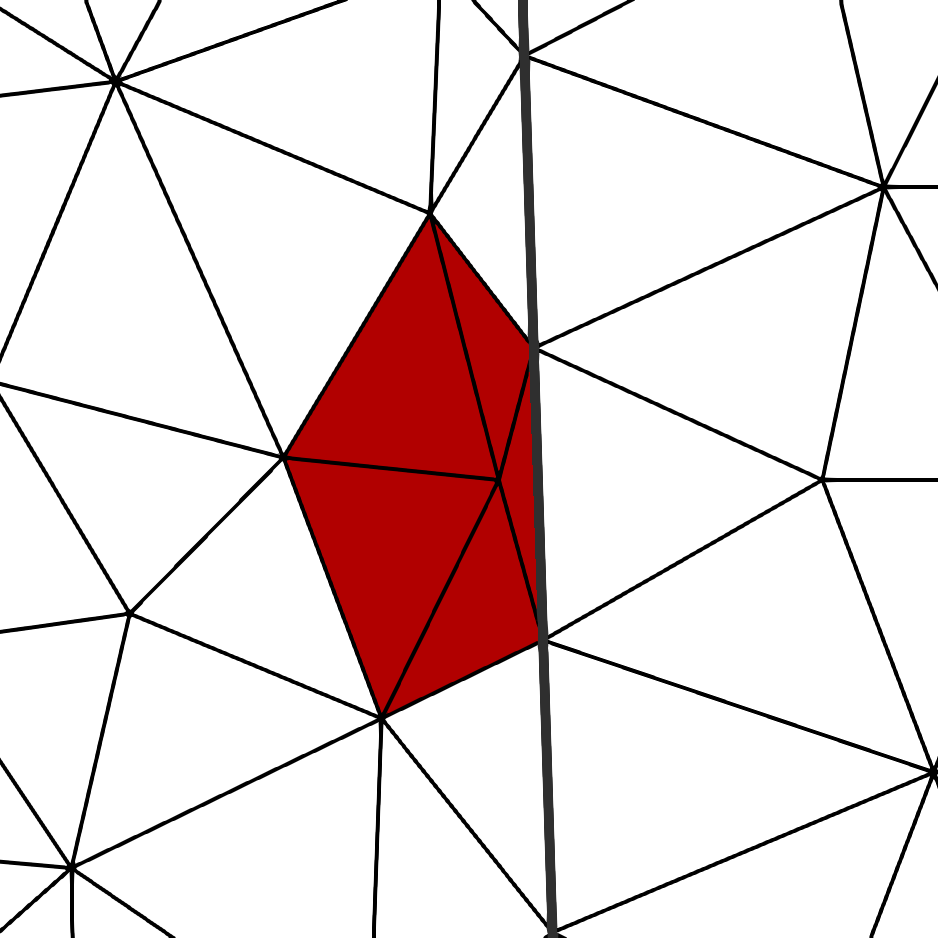}
	\hspace{1mm}
	\includegraphics[width=0.18\textwidth]{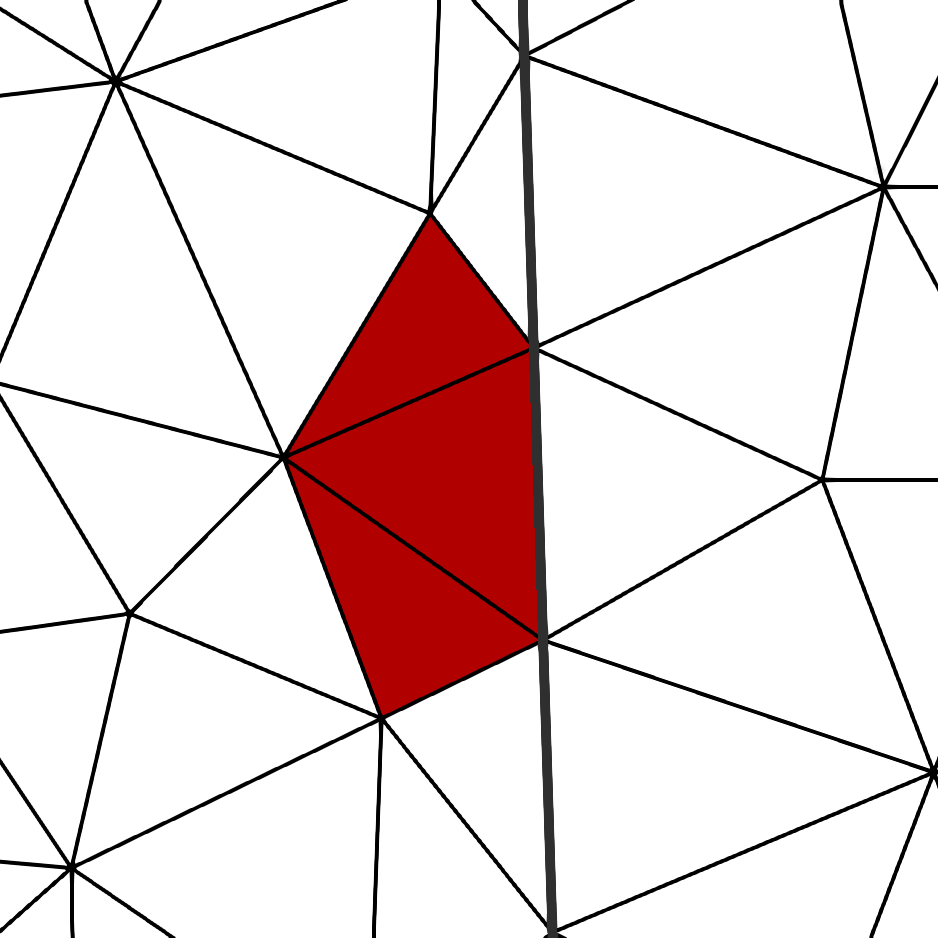}
	\caption{The adaptation of the grid uses connected components for a conservative projection.}
    \label{fig:conncomm}
\end{wrapfigure}

A default adaptation strategy is provided that adapts the triangulation in preparation of vertex movement. It uses an indicator that is defined over several criteria as edge length ratio, radius ratio, edge length and distance to the interface. Per default, re-meshing is performed by retriangulation of holes after removal of vertices and the bisection of edges.
For conservative projection of discrete functions, connected components (compare Fig. \ref{fig:conncomm}) of cells of the old triangulation are constructed that cover the same area as a set of cells in the new one. The weighted average using the exact intersection volumes of cells can be used for conservative projection of cell-wise defined values.
If some edges that belong to the interface have to be refined or coarsened, a similar concept defined on the lower-dimensional triangulation is used.

\section{Numerical Experiments}
\label{sec:experiments}

We demonstrate the performance of the FVMM method in some showcases. For the sake of model validation, we investigate the error between a solution of the reduced model (\ref{eq:ffflowdim})-(\ref{eq:ffffraclowdim}) and a solution of the full dimensional model (\ref{eq:fffgequidim}). In fact, the same scheme can be used to obtain both solutions. In the full dimensional case, the re-meshing capability of \dunemmesh\ is used to track the boundary between bulk and fracture domain where the physical quantities vary.

We consider three numerical experiments for $n=2$. First, we investigate how the FVMM method performs for constant initial data. Second, a prolongating fracture is investigated with source term within the fracture. The solution of the reduced model is compared to the reference solution with the resolved full dimensional fracture. Third, we perform a similar analysis with a squeezing fracture. Finally, we consider a static fracture network for $n=3$.

The source code that was used to produce the results and the raw data of the simulation results is made accessible via DaRUS \cite{darus}.

\subsection{Geometrical Setting and Model Parameters}

We choose a similar geometrical setting for all three cases. Let $\Omega = (0,1)^2$ and $T=1$. The time-dependent fracture $\Gamma(t)$ is given by an ellipse that prolongates and squeezes over time, i.e.,
\begin{align}
  \label{eq:ellipse}
  \Gamma(t) = \{ x_1 = x_2 \mid r(x_1,x_2) \leq R(t) \}, \quad \mathbf{x} = (x_1,x_2)^T.
\end{align}
In (\ref{eq:ellipse}) we have $r(x_1,x_2) \coloneqq \sqrt{ (x_1-0.5)^2+(x_2-0.5)^2 }$ and $R(t) = 0.25 + t v_\text{prolong}$. The aperture is given by
\begin{align}
 d(x_1,x_2,t) = (d_0 - t v_\text{squeeze}) \sqrt{1 - \big(r(x_1,x_2) - R(t)\big)^2}.
\end{align}
The constants $d_0, v_\text{prolong}, v_\text{squeeze} \in \mathbb{R}$ are chosen depending on the case. The geometrical setting is visualized in Fig.\ \ref{fig:2dsetting}. The full dimensional fracture domain $\Omega_f$ is given by $\Omega_f(t) = \{ \mathbf{x} \in \Omega \mid \| \mathbf{x} - \mathbf{s} \| \leq d(\mathbf{x}, t), \mathbf{s} \in \Gamma(t) \}$.

The main parameters are chosen for all three cases as in Table \ref{table:parameters}, if not stated different explicitly. The choice of the parameters is motivated by properties of realistic quantities, but still should be considered as academic. We use $\mathbf{K}_f = \frac{d^2}{12} \mathbf{I}$ inspired by the plane Poisseuille flow.

\begin{figure}[htbp]
\centering
\begin{minipage}{0.5\textwidth}
 \centering
  \begin{tikzpicture}[scale=0.3]
		
  \draw[thick] (0,10) -- (0,0) -- (10,0) -- (10,10);
  \draw[thick] (0,10) -- (10,10);
  \draw[rotate around={45:(5,5)}, gray, fill=gray] (5,5) ellipse (2.5 and 0.5);
  \node at (5,-1){$t=0$};

  \draw[thick] (13,10) -- (13,0) -- (23,0) -- (23,10);
  \draw[thick] (13,10) -- (23,10);
  \draw[rotate around={45:(18,5)}, gray, fill=gray] (18,5) ellipse (4.5 and 0.25);
  \node at (18,-1){$t=1$};

 \end{tikzpicture}
 \captionof{figure}{Geometrical setting where the fracture prolongates and squeezes over time.}
 \label{fig:2dsetting}
\end{minipage}\hfill
\begin{minipage}{0.4\textwidth}
 \centering
 \captionof{table}{Model parameters for all test cases.}
 \begin{tabular}{c|c}
Parameter & Value \\
\hline
$\rho_w$ & $\SI{1000}{\kilo\gram\per\cubic\meter}$ \\
$\rho_{nw}$ & $\SI{500}{\kilo\gram\per\cubic\meter}$ \\
$\mathbf{g}$ & $(0, \SI{-9.81}{\meter\per\second\squared})$ \\
$\mu_w$ & $\SI{1}{\pascal\second}$ \\
$\mu_{nw}$ & $\SI{10}{\pascal\second}$ \\
$\mathbf{K}$ & $\SI{1e-8}{\meter\squared} \ \mathbf{I}$ \\
$\mathbf{K}_f$ & $\frac{d^2}{12} \ \mathbf{I}$ \\
$\phi$, $\phi_f$ & $1$ \\
$k_w(S)$ & $S^2$ \\
$k_{nw}(S)$ & $(1-S)^2$
 \end{tabular}
 \label{table:parameters}
\end{minipage}
\end{figure}

\subsection{Case 1: Constant Initial Data}
\label{subsec:case1}

As first simple benchmark problem we consider a setup where the initial saturations are chosen to be constant $S(\cdot, 0) = S_\Gamma(\cdot, 0) = 1$. The evolution of $S$ is driven by the deformation of the fracture. We apply no external forces ($\mathbf{g} = 0$) and therefore $P \equiv 0$ and $P_\Gamma \equiv 0$. The fracture movement is prescribed by the choices $d_0 = 0.1$, $v_\text{prolong} = 0.25$ and $v_\text{squeeze} = 0$.

\begin{wrapfigure}{r}{0.5\textwidth}
	\includegraphics[width=0.45\textwidth]{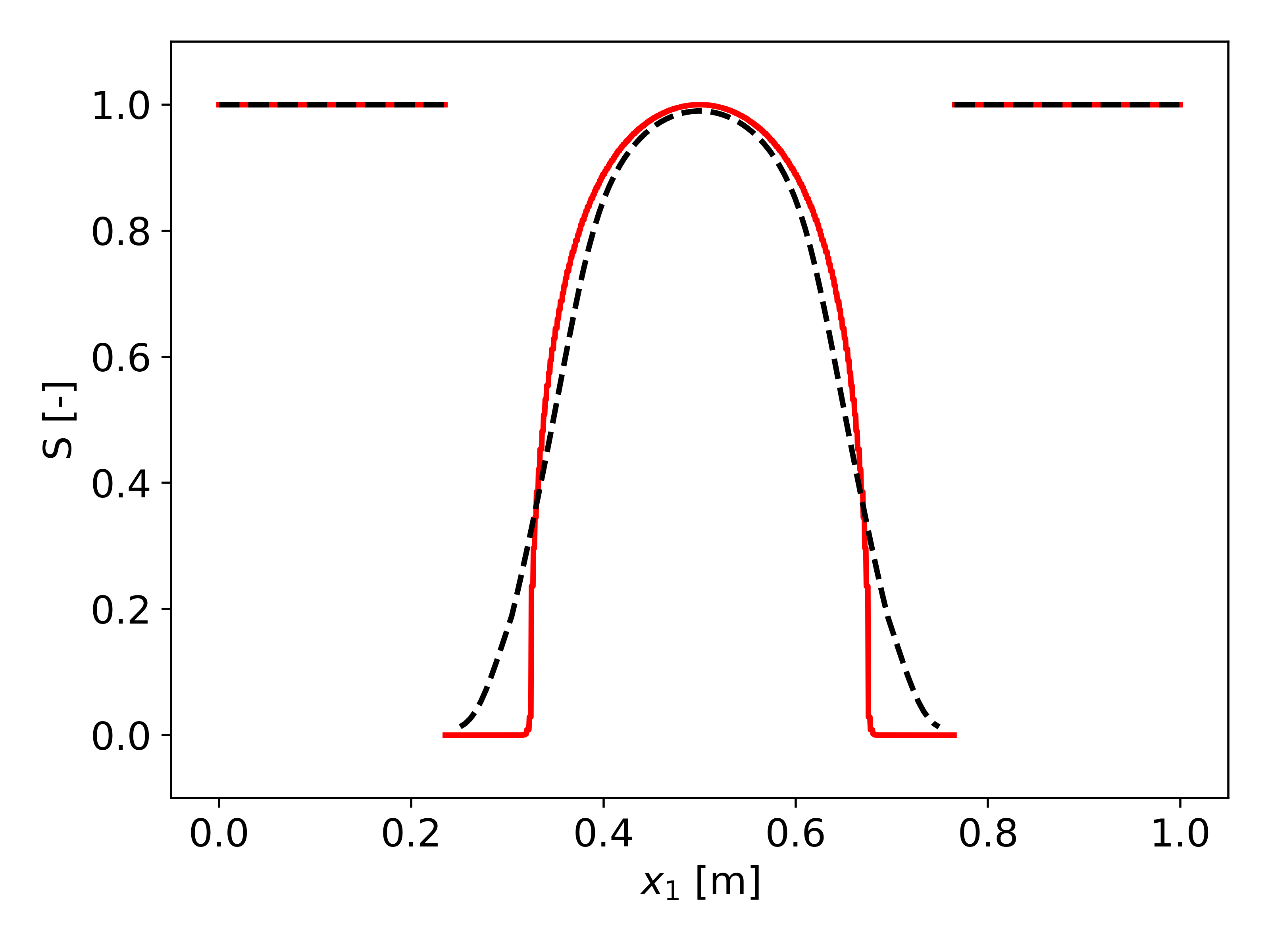}
	\caption{Plot-over-line of saturation for Case 1 at $t = 1.0$. Reduced model (solid) vs. full dimensional model (dashed).}
    \label{fig:case1pol}
\end{wrapfigure}

\begin{figure}[ht]
    \centering
	\includegraphics[width=0.95\textwidth]{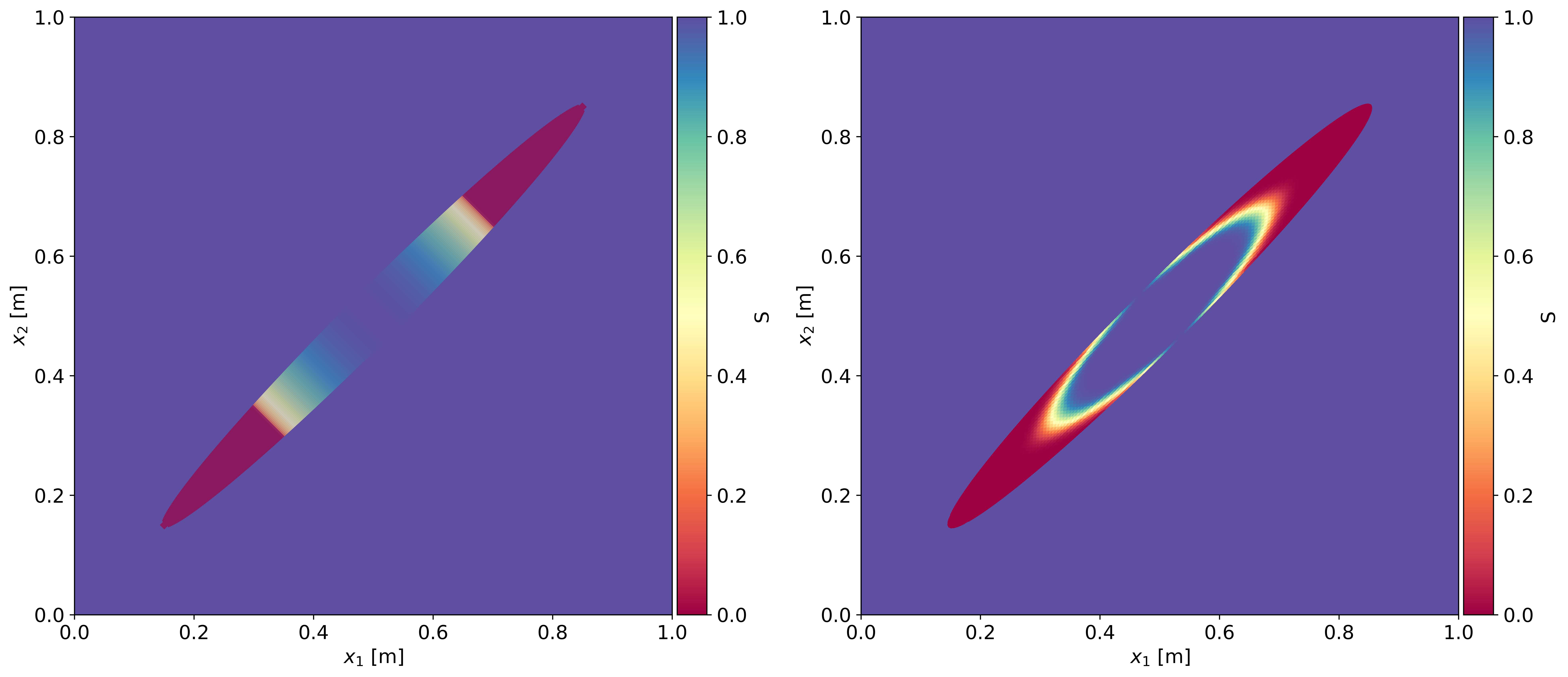}
	\caption{Saturation for Case 1 at $t = 1.0$. Left: Reduced model. Right: Full-dimensional reference.}
    \label{fig:case1}
\end{figure}

A plot-over-line through the fracture center-line is shown in Fig. \ref{fig:case1pol}. Here,  within the fracture, we computed the averaged quantities from the full dimensional reference solution. The result of the saturation at $t = \SI{1}{\second}$ is displayed in Fig. \ref{fig:case1}. In the reduced case, the fracture is visualized as transparent overlay with the corresponding aperture.

We can observe the expected behavior for both the reduced model and the full dimensional model. The saturation drops in direction to the fracture tips. In the full dimensional case, the space occupied by the fracture shows a decreasing saturation over time. This is due to the change of the bulk domain that does not occur in the reduced model.

\subsection{Case 2: A Propagating Fracture}

\begin{wrapfigure}{r}{0.5\textwidth}
    \centering
	\includegraphics[width=0.49\textwidth]{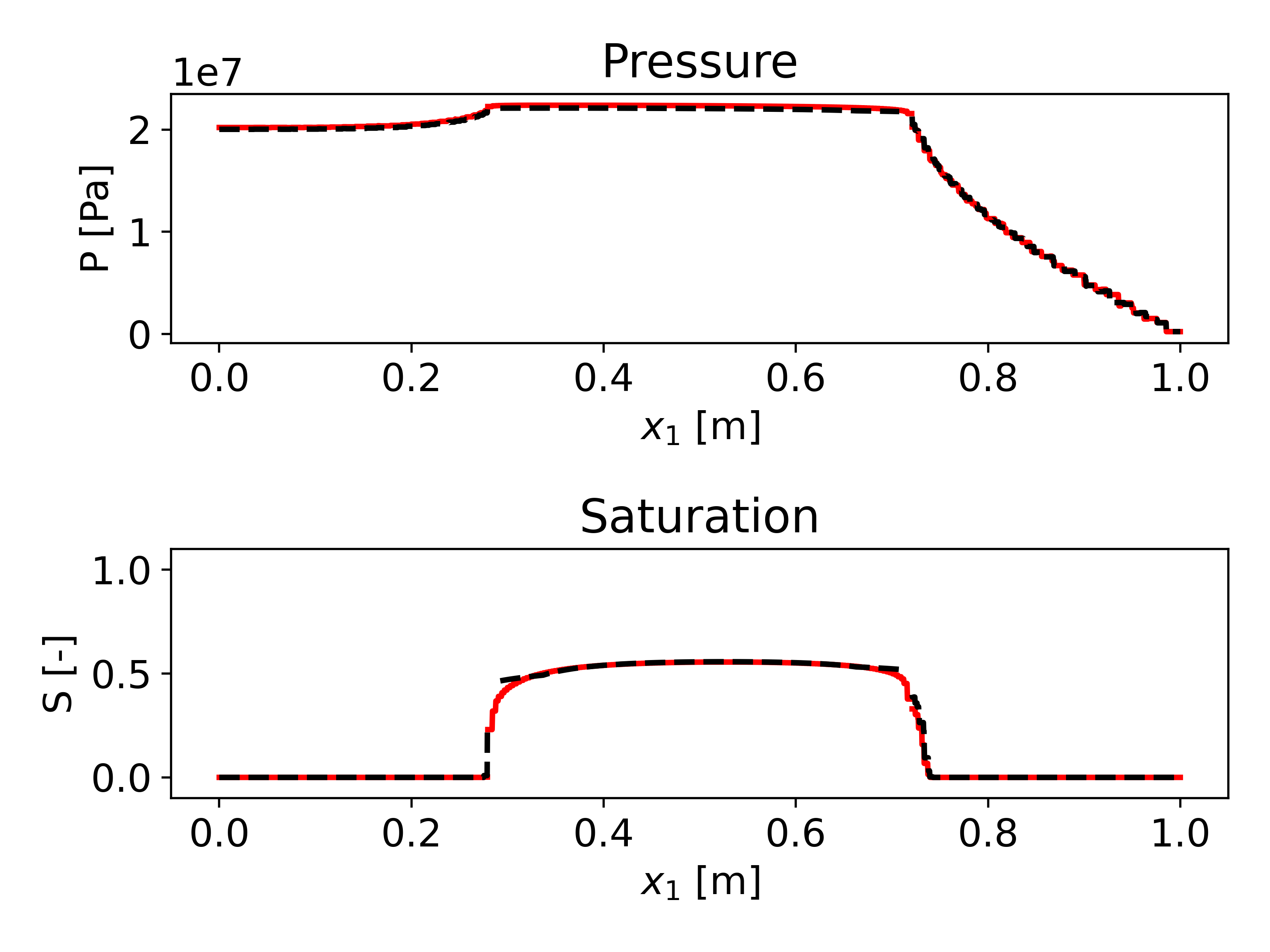}
	\caption{Plot-over-line for Case 2 at $t = 1.0$. Reduced model (solid) vs. full dimensional model (dashed).}
    \label{fig:case2pol}
\end{wrapfigure}

Now, let us consider a propagating fracture with gravity-driven two-phase flow and compare the numerical solution again with a fully-resolved fracture. 

\begin{figure}[ht]
    \centering
	\includegraphics[width=0.95\textwidth]{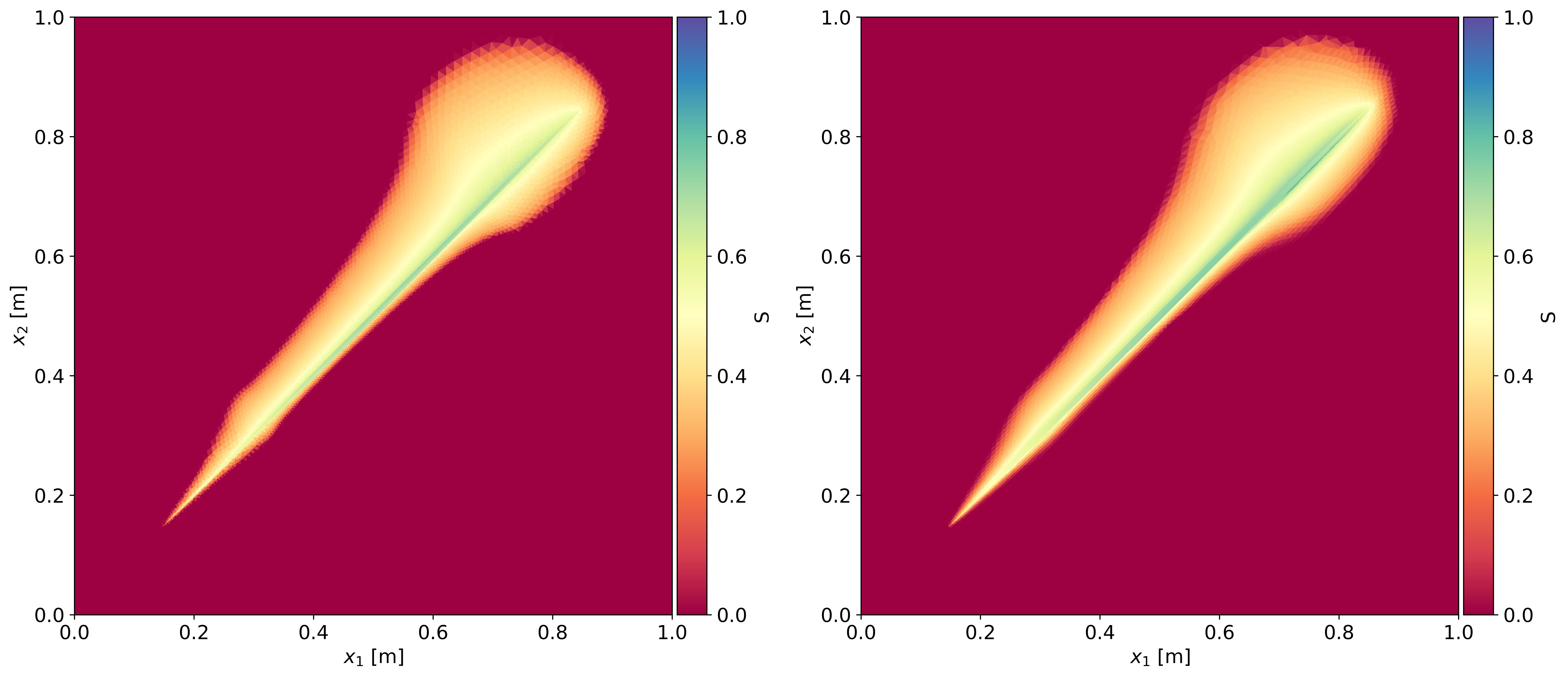}
	\caption{Saturation for Case 2 at $t = 1.0$. Left: Reduced model. Right: Full-dimensional reference.}
    \label{fig:case2}
\end{figure}

The initial saturation is chosen again as $S(\mathbf{x},0)=0$ and $S_\Gamma(\mathbf{s},0) = 0$. The fracture movement is defined by $d_0 = 0.01$, $v_\text{prolong} = 0.25$ and $v_\text{squeeze} = 0$. No-flow boundary conditions are set everywhere except at the top where we fix the pressure to be zero. The source term $q_w^f = q_{nw}^f = \SI{10}{\per\second}$ is applied in $\Gamma$, or $\Omega_f$, respectively.
Again, a plot-over-line through the fracture center-line is displayed in Fig. \ref{fig:case2pol} and the saturation in Fig. \ref{fig:case2}.

We see an overall good agreement of the results between the reduced and the full dimensional model. Small deviations at the fracture tips can be explained by the resolution of the full dimensional grid and the corresponding error in averaging along the orthogonal line segments.

\subsection{Case 3: A Squeezing Fracture}

\begin{wrapfigure}{r}{0.5\textwidth}
    \centering
	\includegraphics[width=0.49\textwidth]{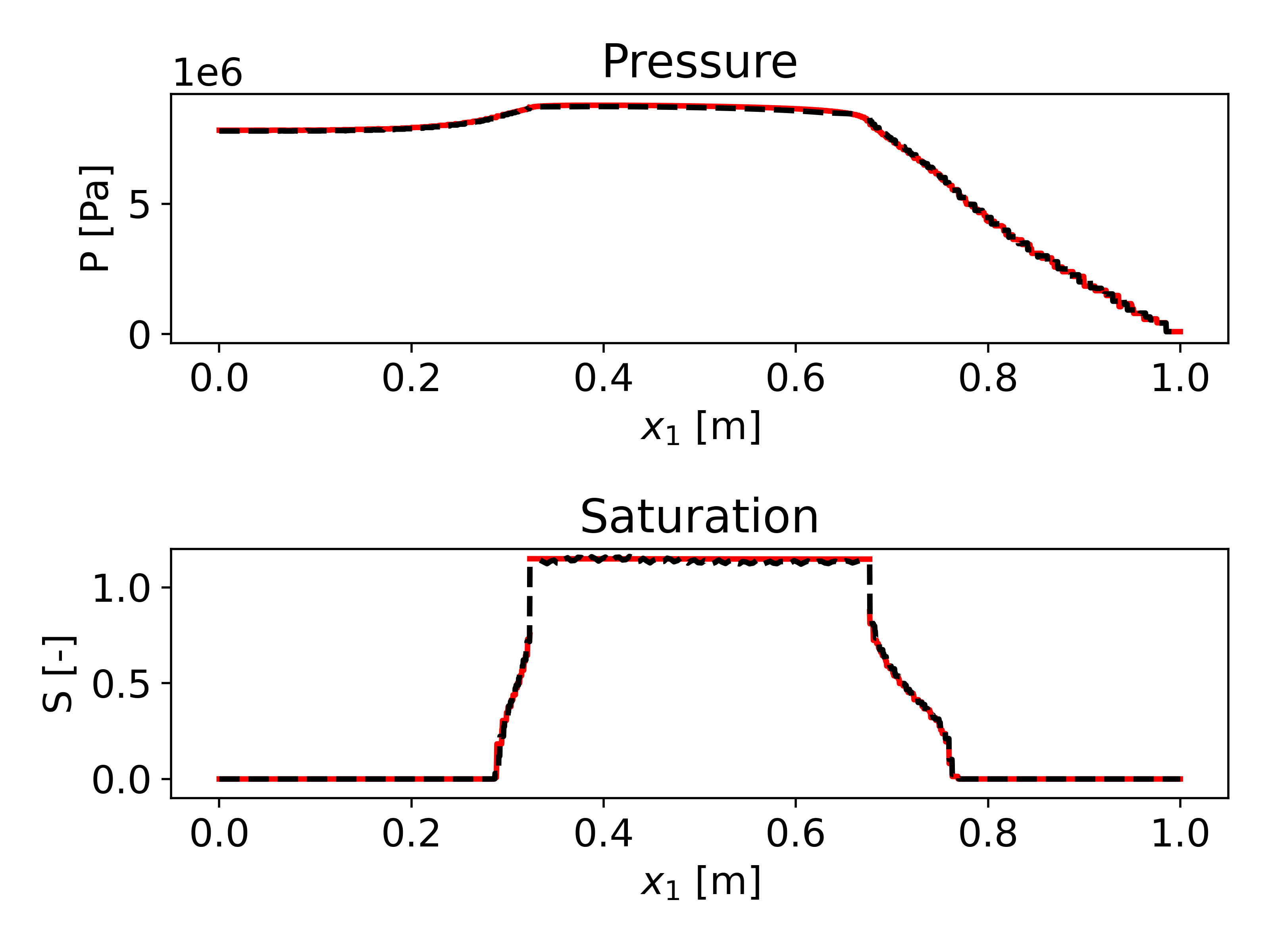}
	\caption{Plot-over-line for Case 3 at $t = 1.0$. Reduced model (solid) vs. full dimensional model (dashed).}
    \label{fig:case3pol}
\end{wrapfigure}

\begin{figure}[ht]
    \centering
	\includegraphics[width=0.95\textwidth]{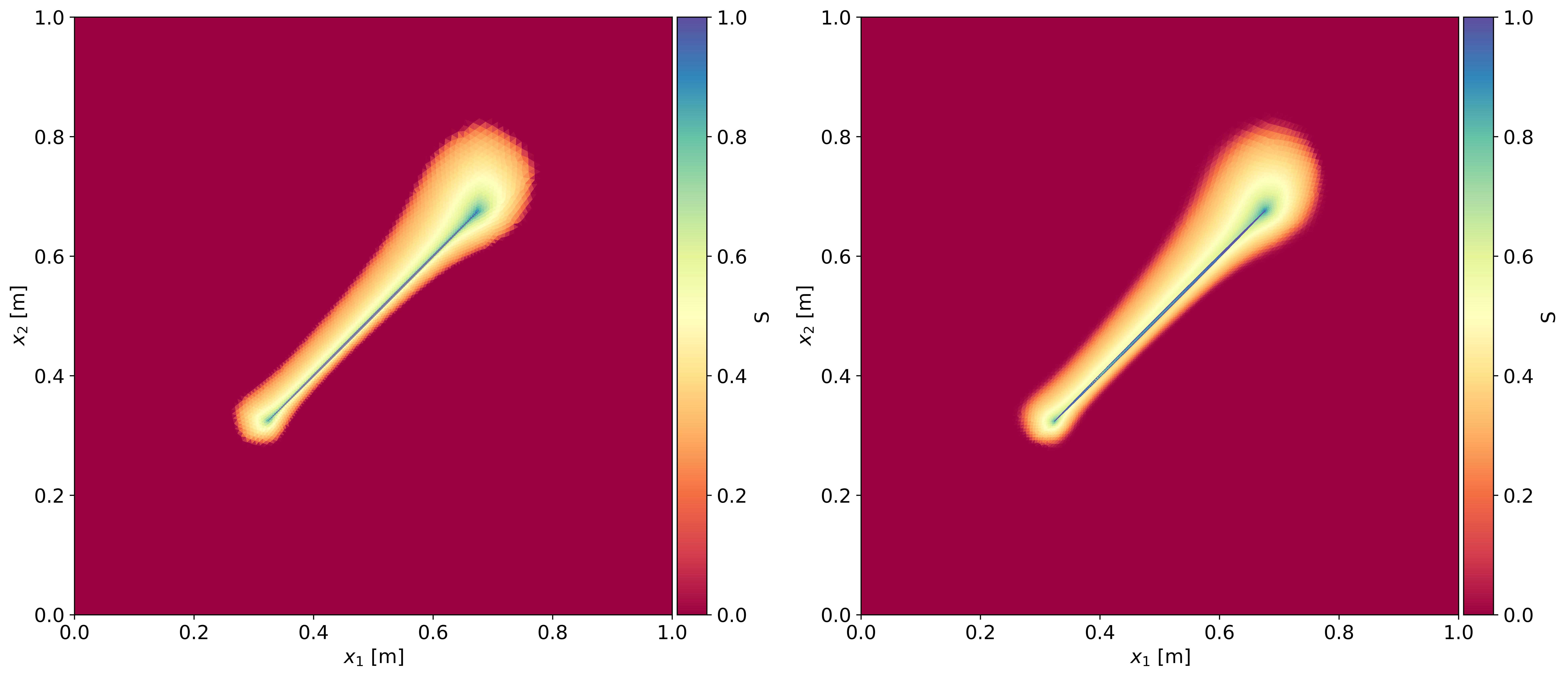}
	\caption{Saturation for Case 3 for $t = 1.0$. Left: Reduced model. Right: Full-dimensional reference.}
    \label{fig:case3}
\end{figure}

Let us consider a squeezing fracture with the physical parameters as in Case 2, but where the fracture movement is defined by $d_0 = 0.01$, $v_\text{prolong} = 0$ and $v_\text{squeeze} = 0.005$.
As in the cases before, we visualize the saturation in Fig. \ref{fig:case3} and a plot-over-line through the fracture center-line in Fig. \ref{fig:case3pol}.

Again, we observe that the solution of the reduced model is in high accordance with the full dimensional reference solution. The saturation profile in the matrix domain around the fractured tips returned in perfect condition.

\subsection{Case 4: Static Fractures in Three Dimensions}

We show a rather exploratory example of 2-dimensional reduced fractures in a 3-dimensional bulk porous medium.
The idea is to show that the proposed reduced model and scheme can also be used in a 3-dimensional setting.
In this example, the fractures are assumed to be static with apertures constant in time.

\begin{figure}[ht]
    \centering
	\includegraphics[width=0.32\textwidth]{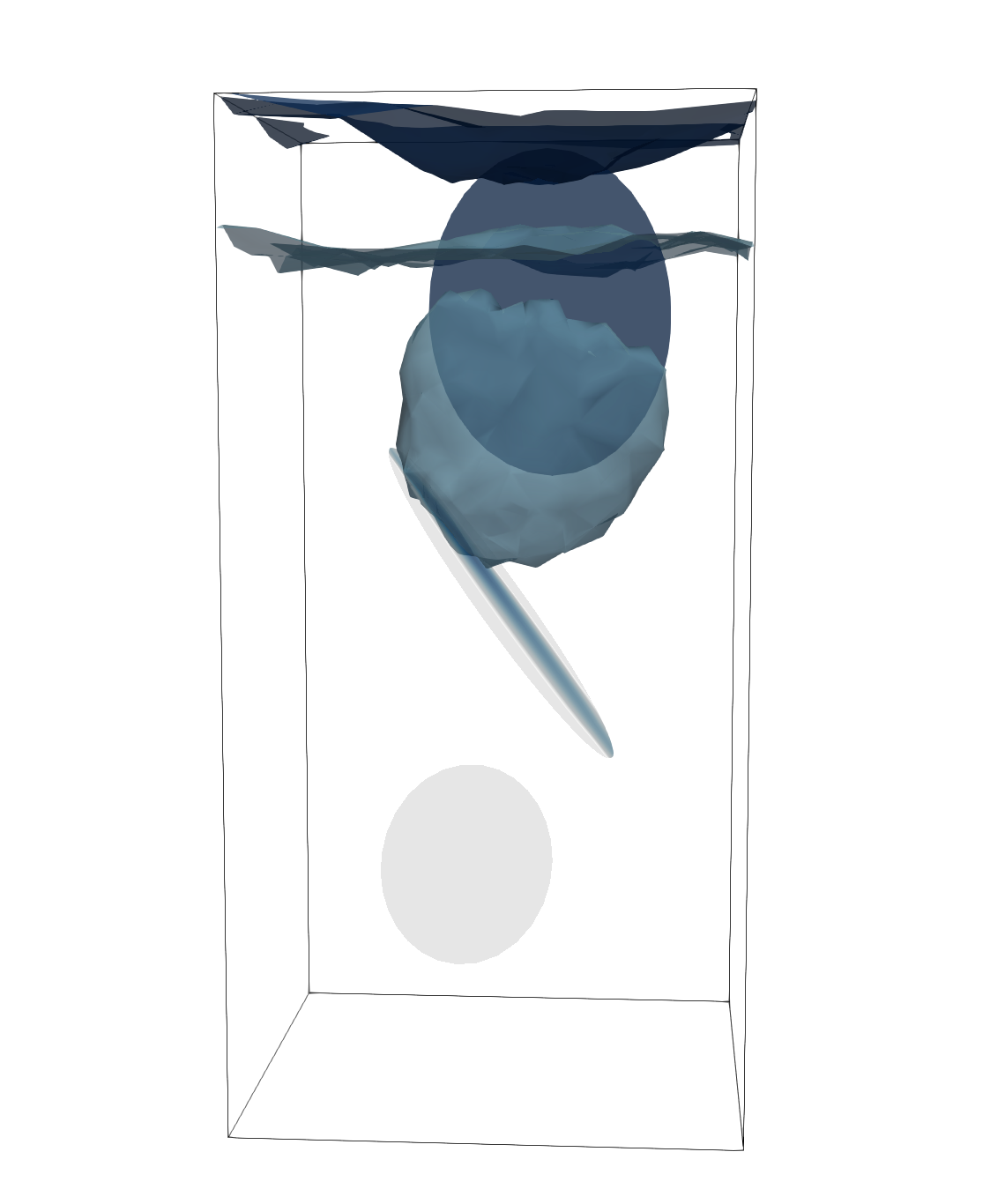}
	\includegraphics[width=0.32\textwidth]{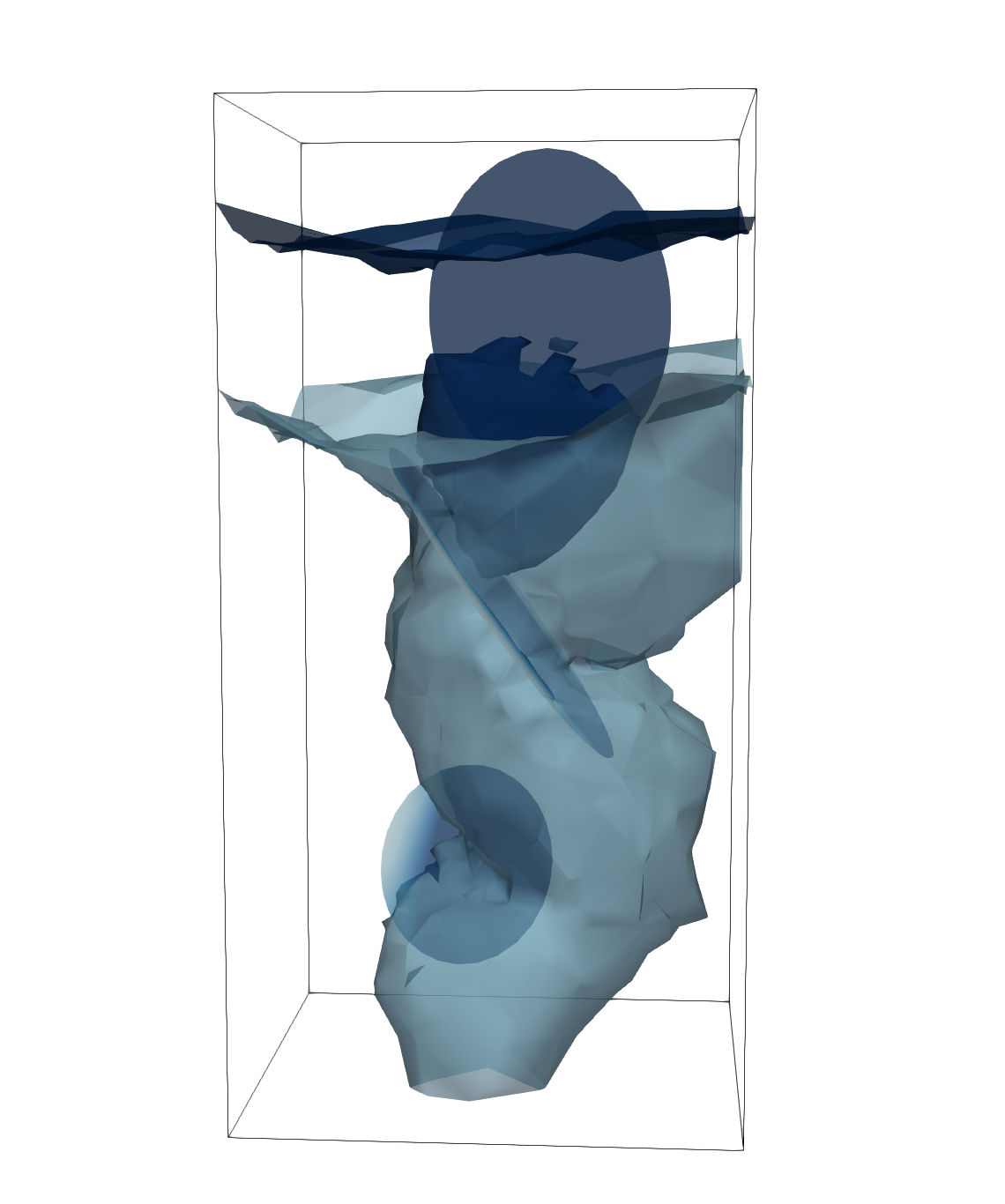}
	\includegraphics[width=0.32\textwidth]{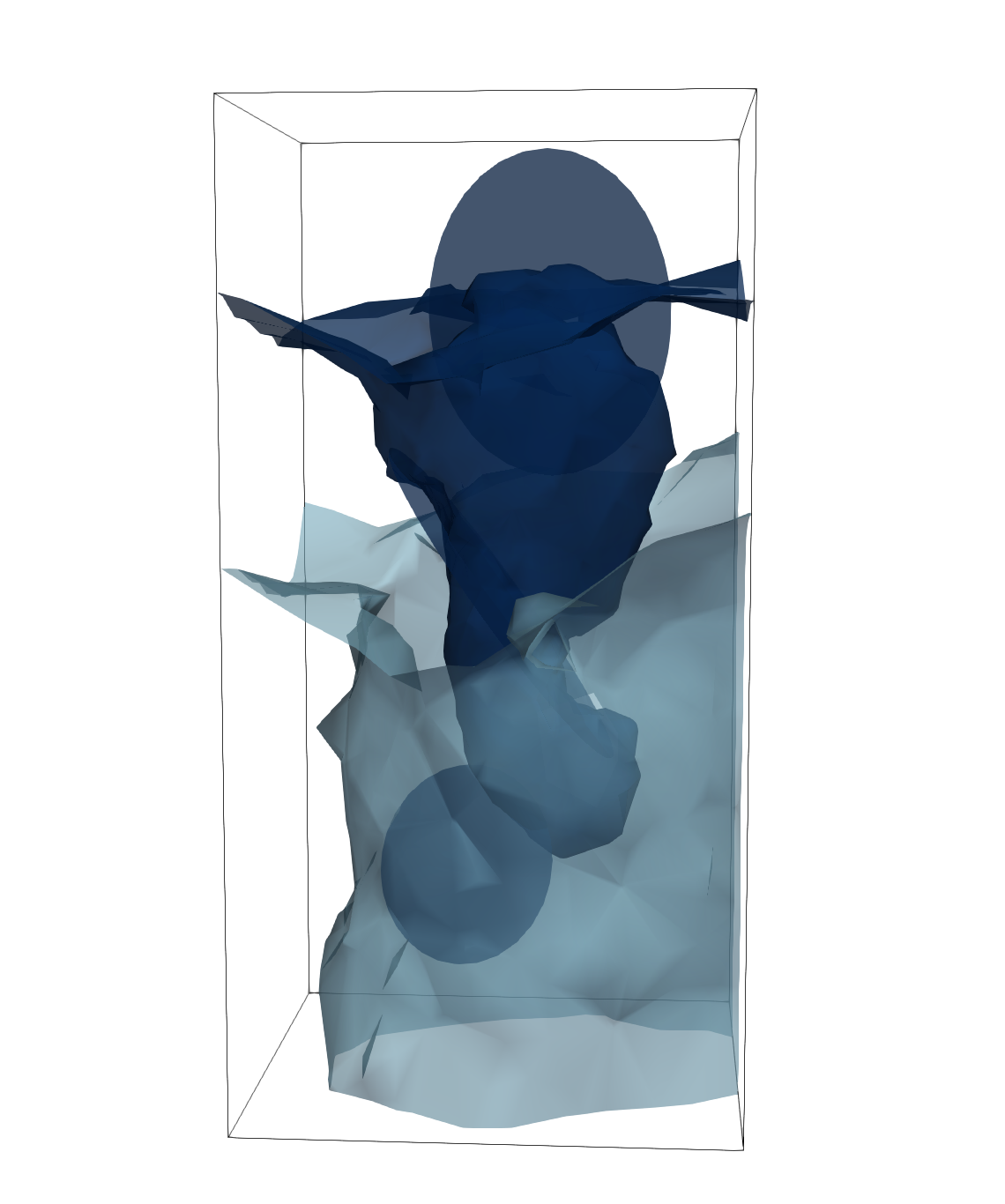}
	\caption{Saturation for Case 4 with 2-dimensional fractures in a 3-dimensional porous medium. Visualized are two contour levels at $S=0.1$ and $S=0.4$ for $t \in \{0.2, 0.6, 1.0\}.$}
    \label{fig:case4}
\end{figure}

\section{Discussion and Outlook}

We presented the derivation of a reduced model for two-phase flow in fractured porous media that takes into account time-dependent fractures and apertures. A numerical method for dynamic fracture propagation scenarios was presented that provides an explicit geometrical representation of the fracture geometries. The method which is able to track movement of lower-dimensional fractures as well as full dimensional fractures was applied different benchmark cases. The results of the reduced and the full dimensional setups were compared and showed good agreement which indicates the validity of the reduced model and the scheme.

Further investigations will aim at more complex geometrical setups with intersecting fractures and topology changes like bifurcations, crossings, etc. In future work, we will also include other physical effects like poro-elasticity and obtain the movement of the fracture tip by the integration of a phase-field model on the microscale locally around the fracture tips \cite{Giovanardi2017, Miehe2015}.

We are working on the generalization of the method to a Discontinuous Galerkin discretization. Work in this direction has been done in e.g. \cite{Antonietti2019}.

\bibliographystyle{unsrtnat}
\bibliography{refs}

\end{document}